\documentclass{article}
\pdfoutput=1

\usepackage{arxiv}

\usepackage[utf8]{inputenc} 
\usepackage[T1]{fontenc}    
\usepackage{hyperref}       
\usepackage{url}            
\usepackage{booktabs}       
\usepackage{amsfonts}       
\usepackage{nicefrac}       
\usepackage{microtype}      
\usepackage{amsmath, amssymb}
\usepackage{graphicx}
\usepackage{latexsym}
\usepackage{color}
\usepackage{tikz}
\graphicspath{{/Users/tranlm/Dropbox/studies/lrecVariance/inst/results/}}
\newcommand\independent{\protect\mathpalette{\protect\independenT}{\perp}}
\def\independenT#1#2{\mathrel{\rlap{$#1#2$}\mkern2mu{#1#2}}}
\definecolor{darkblue}{rgb}{0,0.4,0.9}
\definecolor{gray10}{rgb}{0.1,0.1,0.1}
\definecolor{gray20}{rgb}{0.2,0.2,0.2}
\definecolor{gray30}{rgb}{0.3,0.3,0.3}
\definecolor{gray40}{rgb}{0.4,0.4,0.4}
\definecolor{gray60}{rgb}{0.6,0.6,0.6}
\definecolor{gray80}{rgb}{0.8,0.8,0.8}
\definecolor{gray90}{rgb}{0.9,0.9,.9}
\definecolor{gray95}{rgb}{0.95,0.95,.95}
\definecolor{gray96}{rgb}{0.96,0.96,.96}
\definecolor{lgreen} {RGB}{180,210,100}
\definecolor{dblue}  {RGB}{20,66,129}
\definecolor{ddblue} {RGB}{11,36,69}
\definecolor{lred}   {RGB}{220,0,0}
\definecolor{nred}   {RGB}{224,0,0}
\definecolor{norange}{RGB}{230,120,20}
\definecolor{nyellow}{RGB}{255,221,0}
\definecolor{ngreen} {RGB}{98,158,31}
\definecolor{dgreen} {RGB}{78,138,21}
\definecolor{nblue}  {RGB}{28,130,185}
\definecolor{jblue}  {RGB}{20,50,100}
\definecolor{nnyellow}{RGB}{235,200,0}
\definecolor{purple}{RGB}{150, 0, 120}
\definecolor{sgGreen} {RGB}{20, 180, 50}
\definecolor{revised}{rgb}{0,0,0.9}

\newcommand{\openr}{\hbox{${\rm I\kern-.2em R}$}}
\newcommand{\openn}{\hbox{${\rm I\kern-.2em N}$}}

\newcommand\numberthis{\addtocounter{equation}{1}\tag{\theequation}}

\title{Robust variance estimation and inference for causal effect estimation}
\author{
  Linh Tran,
  \quad
  Maya Petersen,
  \quad
  Joshua Schwab,
  \quad
  Mark J van der Laan \\
  School of Public Health \\
  University of California, Berkeley \\
  Berkeley, CA 94720 \\
  \texttt{tranlm@google.com}, \texttt{mayaliv@berkeley.edu},
  \texttt{joshuaschwab@yahoo.com}, \texttt{laan@berkeley.edu} \\
}

\begin{document}
\maketitle

\begin{abstract}
We consider a longitudinal data structure consisting of baseline covariates,
time-varying treatment variables, intermediate time-dependent covariates, and a
possibly time dependent outcome. Previous studies have shown that estimating the
variance of asymptotically linear estimators using empirical influence functions
in this setting result in anti-conservative estimates with increasing magnitudes
of positivity violations, leading to poor coverage and uncontrolled Type I
errors. In this paper, we present two alternative approaches of estimating the
variance of these estimators: (i) a robust approach which directly targets the
variance of the influence function as a counterfactual mean outcome, and (ii) a
non-parametric bootstrap based approach that is theoretically valid and lowers
the computational cost, thereby increasing the feasibility in non-parametric
settings using complex machine learning algorithms. The performance of these
approaches are compared to that of the empirical influence function in
simulations across different levels of positivity violations and treatment
effect sizes.
In the appendix, we generalize the robust approach of estimating variance
to working marginal structural working models.
\end{abstract}


\section{Introduction}

A number of estimators are available for the treatment specific mean outcome
parameter (and the corresponding causal contrasts) based on longitudinal data
structures, such as inverse probability weighting (IPW) \citep{Horvitz1952,
Robins1998}, double robust augmented IPW (AIPW) \citep{Robins1992, Robins1994,
Robins2000c, Robins2000, Robins2001, Rotnitzky2005}, and targeted minimum loss-based
estimation (TMLE) \citep{vanderLaan2011a}. Variance estimation for each of
these are conventionally achieved by using their corresponding influence
functions (IF) on the empirical distribution or by resampling methods such as
the non-parametric bootstrap. However, a number of shortcomings exists with
these variance estimation approaches. In particular, no theory for exists for
the non-parametric bootstrap when using data adaptive methods for estimation of
nuisance parameters, and both IF-based and bootstrap based confidence intervals
can become anti-conservative with increasing levels of practical positivity
violations. For example, van der Laan et al. \citep{vanderLaan2011a} found IF-based
variance estimates for the intervention specific mean outcome that were
anti-conservative when compared with the Monte-Carlo variance of the TMLE,
leading to invalid confidence intervals. Petersen el al. \citep{Petersen2012}
found poor coverage for influence function-based confidence intervals, owing to
both a result of practical positivity violations and relatively rare outcomes.
This behaviour is especially true under sparsity in finite samples, even when
the assumptions for asymptotic validity of these estimators hold
\citep{Petersen2012}. As a consequence, statistical inference based on these
estimators of variance becomes unreliable when the treatment mechanism
practically or theoretically violates the underlying positivity assumption.

Additionally, under sparsity issues, the estimated variance may also fail to
raise a red flag for unreliable statistical inference \citep{Petersen2012}. For
example, these estimates of the asymptotic variance are not sensitive to
theoretical violations of the positivity assumptions under which the asymptotic
variance would be infinity, i.e. when positivity fails. Consequently, it is less
likely that the analyst will be able to determine if the data at hand provides
insufficient information to estimate the desired causal parameter with any
reasonable degree of accuracy.

Previous work \citep{Petersen2012, Petersen2014c} proposed estimating the
asymptotic variance of the estimator with a parametric bootstrap-based on a fit
of the density of the data generating distribution, involving estimation of the
treatment mechanism and the $G$-computation factor of the likelihood. This
proposal corresponds with evaluation of the variance of a given estimator using
the data at hand as a given data generating experiment. The consistency of this
estimator relies on correct specification of both the treatment mechanism and
the $G$-computation factor. This parametric bootstrap integrates over sparse
events and therefore will explode the variance. An extremely large number of
samples is therefore needed to get the true variance under this Monte Carlo
scheme. As a consequence, this parametric bootstrap-based variance estimate was
only proposed as a measure to raise a red flag for unreliable statistical
inference. In addition, in the context of sparsity, one needs to sample a large
number of bootstrap samples and refit the likelihood in each iteration in order
to obtain a valid evaluation of the estimator variance, in order to capture the
rare observations that nonetheless heavily contribute to this variance. Thus,
this semi-parametric bootstrap method is extremely computer intensive, making
this Monte Carlo scheme an intractable method for complex estimators and complex
data generating distributions.

In this article we use analytic expressions to compute the variance of the
efficient influence function (EIF) \citep{Hampel1974, Robins1992} which provide
the asymptotic variance of estimators solving the estimating equation
corresponding to this function. These analytic expressions naturally integrate
over the rare observations, and thereby avoid the finite sample bias in variance
estimation using standard influence curve or non-parametric bootstrap based
methods due to rare observations mentioned above. With this, we construct plug-
in type estimators of these asymptotic variances that are consistent if both the
treatment mechanism and treatment specific means of specified outcomes are
consistently estimated. These estimators require estimation of the treatment
mechanism and several treatment specific means of specified outcomes (defined as
a function of the observed data structure, indexed by the estimator of the
treatment mechanism), which can be estimated with either an estimating equation
type IPW estimator or an efficient double robust method such as a targeted
minimum loss-based estimator (TMLE). The resulting variance estimator, unlike
current alternatives based on taking the variance of the empirical influence
function, or using a non-parametric bootstrap, will become very large whenever
the estimated treatment mechanism reflects practical or theoretical violations
of the positivity assumption.

While this newly presented approach performs well in estimating the asymptotic
variance of estimators solving the estimating equation corresponding to the EIF,
a lower finite sample variance should be expected for substitution based
estimators such as TMLE \citep{vanderLaan2011a}, due to the guaranteed parameter
boundaries provided by the estimator. We therefore additionally present a second
bootstrap based approach of estimating the finite sample variance. The approach
does not require re-estimation of the treatment mechanisms and the $q$-factor of
the likelihood and therefore reduces the computational burden. It is
asymptotically consistent under reasonable assumptions; namely, the same
essential assumptions needed for the estimator of the target parameter itself to
be asymptotically linear. The resulting reduction in the computational load
(compared to a fully non-parametric bootstrap approach which refits the
likelihood for each iteration) allows for a more tractable approach at
estimating the variance.

\subsection{Organization of this paper}

In Section \ref{sec:data}, we formally define the observable data, likelihood,
and statistical model for its distribution. Our target parameter of the
treatment specific mean is defined along with its EIF. We briefly review the
causal model and assumptions under which this statistical quantity corresponds
with the desired causal parameter of the counterfactual distribution, along with
the currently common approach of influence function (IF) based estimator
variance estimation.

Section \ref{sec:robustVar} presents an approach for robust estimation of the
variance of the EIF under sparsity. The expression for the variance of the
efficient influence function is presented along with both an IPW and TMLE based
approach at estimating this parameter. To help illustrate, an example is given
for a point treatment setting under a static treatment regime. Advantages of
this new approach are covered. The Appendix generalizes the approach to working
marginal structural working models and provides proofs.

Section \ref{sec:bootstrap} discusses the second approach of estimating the
estimator variance using the bootstrap, using a modified TMLE. This bootstrap
approach does not require re-estimation of the treatment mechanisms and the
$q$-factor of the likelihood, therefore reducing the computational intensity
required.

Section \ref{sec:sim} illustrates the performance of the variance estimators
presented in Sections \ref{sec:robustVar} and \ref{sec:bootstrap} by applying
them in simulations to both a single time-point and longitudinal setting.
Results show that the robust approach at estimating variance is conservative for
substitution based estimators of the mean outcome, while the bootstrap approach
results in estimates close to the observed Monte-Carlo variance. The resulting
confidence intervals are demonstrated to be valid under the newly proposed
variance estimation approaches, while the bootstrap approach is shown to retain
higher statistical power.


We conclude with a discussion in Section \ref{sec:disc}, which reviews the
results, benefits of this new approach, potential limitations, and future
directions.

\section{Definition of data and statistical estimation}
\label{sec:data}

Consider a longitudinal study in which subjects are seen at each time point
$t$ from $t=0,1,\ldots,K+1$. The observable data structure on a randomly sampled
subject is
\begin{equation}
  O=(L(0), A(0), L(1), A(1), \ldots, A(K), Y=L(K+1)) \overset{iid}{\sim} P_0
\end{equation}

where $L(0)$ includes all baseline covariates, $A(t)$ denotes an intervention
node at time $t$, and $L(t)$ denotes all time-varying covariates at time point
$t$, measured between the intervention nodes $A(t^{-})$ and $A(t)$, where for
notational convenience we define $t^{-} \equiv t-1$. Our outcome of interest
$Y=L(K+1)$ is an outcome measured after the final treatment $A(K)$. We observe
$n$ independent and identically distributed (iid) copies copies
$O_i:i=1,\ldots,n$, of $O$.

The likelihood $L(O)$ for the observable data is the product of conditional
probabilities such that the likelihood for subject $i$ is
\begin{equation}
\label{eq:likelihood}
\begin{split}
  L(O_i) &= p_0(L_i(0),A_i(0),L_i(1),A_i(1),\ldots,L_i(K+1)) \\
    & = p_0(L_i(K+1)|\bar{L}_i(K),\bar{A}_i(K))\cdot p_0(A_i(K)|\bar{L}_i(K),\bar{A}_i(K-1)) \\
      &\hspace{.15in}\cdot p_0(L_i(K)|\bar{L}_i(K-1),\bar{A}_i(K-1)) \cdot p_0(A_i(K-1)|\bar{L}_i(K-1),\bar{A}_i(K-2)) \\
      &\hspace{.15in}\cdots p_0(L_i(0)) \\
    & = \left[\prod_{t=0}^{K+1}
    \underbrace{p_0(L_i(t)|\bar{L}_i(t^{-}),\bar{A}_i(t^{-}))}_{q_{0,t}(L(t)|Pa(L(t))}\right]
    \cdot \left[\prod_{t=0}^{K}
    \underbrace{p_0(A_i(t)|\bar{L}_i(t),\bar{A}_i(t^{-}))}_{g_{0,t}(A(t)|Pa(A(t))}\right]
\end{split}
\end{equation}

where $\bar{X}(t) \equiv (X(1),X(2),\ldots,X(t))$,
$A(-1)=L(-1)=\varnothing$, and $p_0(o)$ denotes $p_0(O=o)$ under the true
distribution $P_0$ where we assume $O$ is discrete for sake of presentation.

The statistical model ${\cal M}$ for the data involves assumptions, if any, only
on the conditional distributions of $A(t)$, given
$Pa(A(t))=(\bar{L}(t),\bar{A}(t^{-}))$, $t=0,\ldots,K$. Let
\[P^{d}_{0}(l)\equiv\prod_{t=0}^{K+1} P_{0,L(t)}(l(t)\mid
\bar{l}(t^{-}),d(\bar{l}(t^{-}))\] denote the $G$-computation formula for the
post-intervention distribution of an intervention that sets
$\bar{A}(K)=d(\bar{l}(K))$ \citep{Robins1986}. We use the notation $P_{L(t)}$
for a conditional distribution of $L(t)$, given
$Pa(L(t))=(\bar{L}(t^{-}),\bar{A}(t^{-}))$. Let
$L^{d}=(L(0),\ldots,Y^{d}=L^{d}(K+1))$ be a random variable under the
post-intervention distribution $P^{d}_0$. The statistical target estimand is
defined here as $\Psi(P_0)=\mathbb{E}_{P^{d}_0}[Y^{d}]$, i.e. the mean of the
outcome at time $K+1$ under this distribution. We note that $\Psi:{\cal
M}\rightarrow\openr$ represents a target parameter mapping on the statistical
model to the real line. Defining $t^{+}\equiv t+1$, the EIF of $\Psi$ at $P$ is
given by \[D^*(P)(O)=\sum_{t=0}^{K+1}D_t^*(P)(O),\] \citep{Robins2000,Bang2005}
where

\begin{eqnarray*}
D_0^*(P)(L(0)) &=&\bar{Q}^{d}_1-\bar{Q}^{d}_0\\
D_t^*(P)(\bar{A}(t^{-}),\bar{L}(t^{-}))&=&
H_t(g)(\bar{Q}^{d}_{t^{+}}-\bar{Q}^{d}_t):\mbox{$t=1,2,\ldots,K+1$}\\
\end{eqnarray*}
where
\begin{eqnarray}
\label{eq:H}
H_t(g)(\bar{A}(t^{-}),\bar{L}(t^{-}))&=&\frac{\mathbb{I}(\bar{A}(t^{-})=d(\bar{l}(t^{-})))}{g_{0:t^{-}}(\bar{A}(t^{-}),\bar{L}(t^{-}))}\\
\bar{Q}^{d}_{K+2}&=&Y \\
\bar{Q}^{d}_t(\bar{L}(t^{-}))&=&\mathbb{E}_{P}[Y^{d}\mid
\bar{L}^{d}(t^{-})=\bar{L}(t^{-})]:\mbox{$t=1,2,\ldots,K+1$}\\
\bar{Q}^{d}_0&=&\mathbb{E}_P[Y^{d}].
\end{eqnarray}

It should be noted that $g_{0:t^{-}}(\bar{A}(t^{-}),\bar{L}(t^{-}))$ represents
the cumulative probability of treatment up to time $t-1$ and that
$\bar{Q}^{d}_t(\bar{L}(t^{-}))=\mathbb{E}_P[\bar{Q}^{d}_{t^{+}}\mid
\bar{L}(t^{-}),\bar{A}(t^{-})=d(\bar{l}(t^{-}))]$ is defined by recursive regression,
starting at $t=K+1$ and moving backwards in time. For notational convenience, we
let $H_0=1$ so that
\[D^*(P)(O)=\sum_{t=0}^{K+1}H_t(g)(\bar{Q}^{d}_{t^{+}}-\bar{Q}^{d}_t).\]

\subsection{Causal model}

Under additional assumptions about our data generating where our target
statistical estimand is equal to the mean of the counterfactual outcome $Y_{d}$
under intervention to set the vector of treatment nodes to value
$d(\bar{l}(K))$. Specifically, for interventions of interest $d \in {\cal D}$,
we assume sequential randomization \citep{Robins1986} \[ Y_{d} \independent A(t)
| \bar{L}(t),\bar{A}(t^{-}): t=0,1,..,K \] and positivity
\citep{Robins1999b} \[
P(\bar{A}(t)=d(\bar{l}(t))|\bar{L}(t),\bar{A}(t^{-})=d(\bar{l}(t)))>0
\text{ a.e.}:t=0,1,\ldots,K \] Regarding the assumption of positivity, we note
that as $P(\bar{A}(t)=d(\bar{l}(t))|\bar{L}(t),\bar{A}(t^{-})=d(\bar{l}(t)))\rightarrow
0$, we have that $H_t(g) \rightarrow \infty$ resulting in $\text{var}[D^*(P)(O)]
\rightarrow \infty$ and the previously mentioned practical and theoretical
positivity violations.

\subsection{Review of Influence Function based variance}

Recall that an estimator $\hat{\Psi}(P_n)$ is considered to be asymptotically
linear if and only if \[ \hat{\Psi}(P_n) - \Psi(P_0) =
\frac{1}{n}\sum_{i=1}^nD(P_0)(O_i)+o_p(n^{-1/2})\] for some mean $0$ finite
variance influence function $D(P_0)(O)$ \citep{Hampel1974}. If an estimator is
asymptotically linear, then it will be asymptotically normal with variance equal
to the variance of the influence function over $n$. The asymptotic variance of
the estimator can therefore be consistently estimated with the variance of the
empirical influence function $D(P_n)(O)$, i.e.
$\hat{\text{var}}[\hat{\Psi}(P_n)]=\text{var}[D(P_n)(O)]/n$, which implies an
asymptotically valid confidence interval.

\subsubsection{Targeted minimum loss based estimation (TMLE)}

One such estimator that solves the the estimating equation corresponding to the
efficient influence function for intervention specific mean outcomes is Targeted
Minimum Loss-based Estimation \citep{vanderLaan2011a}. This estimator
solves the estimating equation by forming an intial fit of the $q_0$ portion of
the likelihood and subsequently perturbing it such that the estimating equation
is solved. We assume use of this estimator for our estimation problem, such that
our attention is focused on estimation of the estimator's variance. Note that
our proposed variance estimators also apply to estimating equation approaches,
such as the double robust augmented IPW (AIPW) \citep{Robins1992, Robins1994,
Robins2000c, Robins2000, Robins2001, Rotnitzky2005}.


\section{Semi-targeted estimation of the EIF variance}
\label{sec:robustVar}

We can directly target the variance of $D(P_{0})(O)$ as an expectation, allowing
us to estimate the variance as a mean. The following describes how to obtain a
TMLE of the variance of each component of the EIF $\sigma_{t}^2$ in the setting
of a scalar parameter. We provide a proof for the more general working MSM
setting in the Appendix for the interested reader.

\subsection{Expression for variance of the EIF for  \texorpdfstring{$\mathbb{E}Y_{d}$}{Lg}}

Under regimens $d(\bar{l}(K))$, we have
\[
\sigma^2_0 \equiv \mathbb{E}_0[D^*(P_0)(O)]^2=\sum_{t=0}^{K+1} \mathbb{E}_0
[H_t^2(g_0) (\bar{Q}^{d}_{0,t^{+}}-\bar{Q}^{d}_{0,t})^2].
\]
Using the expression for $H_t(g)$ from Equation \eqref{eq:H}, and first taking
the conditional expectation w.r.t. $\bar{A}(t^{-})$ given
$X=(\bar{L}^{d}:d)$, it follows that this can be written as:
\begin{equation}
\label{eq:eicvar1}
\sigma^2_0 = \sum_{t=0}^{K+1} \mathbb{E}_{P_0^{d}}
\left[\frac{(\bar{Q}^{d}_{0,t^{+}}-\bar{Q}^{d}_{0,t})^2
(\bar{L}^{d}(t))}{g_0(d(\bar{l}(t^{-})),\bar{L}^{d}(t^{-}))}\right],
\end{equation}

where we define $g_0(d(\bar{l}(-1)),\bar{L}^d(-1))=1$ so that the term at $t=0$
equals
$\mathbb{E}_{L(0)}[\bar{Q}^{d}_{0,1}(L(0))-\mathbb{E}_0Y^{d}]^2$.
This is simply a sum of expectations over $t \in \{0,1,\ldots,K+1\}$. For
notational convenience, we re-write Equation \eqref{eq:eicvar1} as
\begin{equation}
\label{eq:eifvar2}
\sigma^2_0=\sum_{t=0}^{K+1}\sigma^{2,d}_t = \sum_{t=0}^{K+1}
\mathbb{E}_{P_0^{d}}\left[S_t^{d}(\bar{Q}_0,g_0)(\bar{L}^{d}(t))\right]
\end{equation}
for the specified function
\[ S_t^{d}(\bar{Q}_0,g_0)(\bar{L}^{d}(t))
\equiv
\frac{(\bar{Q}^{d}_{0,t^{+}}-\bar{Q}^{d}_{0,t})^2(\bar{L}^{d}(t))}
{g_0(d(\bar{l}(t^{-})),\bar{L}^{d}(t^{-}))}:
\mbox{$t=0,1,\ldots,K+1$}
\]

Note that, given ($\bar{Q}_0,g_0$), $\mathbb{E}_{P_0^{d}}
S_t^{d}(\bar{Q}_0,g_0)$ is the mean of a counterfactual
$S^{d}_t(\bar{Q}_0,g_0)(\bar{L}^{d}(t))$, i.e., the mean of a real valued
function (indexed by $d(\bar{l})$ itself) of $\bar{L}^{d}(j)$, which needs to be
estimated based on the longitudinal data structure $L(0),A(0),\ldots,
A(t-1),L(t)$. Given $\bar{Q}_0,g_0$, we observe the outcome
$S^{d}_t(\bar{Q}_0,g_0)(\bar{L}_i(t))$, $i=1,2,\ldots,n$, so that we can
represent the observed data structure as
$L(0),A(0),\ldots,A(t-1),S^{d}_t(\bar{Q}_0,g_0)(\bar{L}(t))$, and we wish to
estimate the statistical target parameter
\begin{equation}
\label{eq:varparm}
\mathbb{E}_{P_0^{d}}
S^{d}_t(\bar{Q}_0,g_0)=\sum_{\bar{l}(t)}
S^{d}_t(\bar{Q}_0,g_0)(\bar{l}(t))P_0^{d}(\bar{L}^{d}(t)=\bar{l}(t)):t=0,1,\ldots,K+1,
\end{equation}
where again we assume $l(t)$ is discrete for sake of presentation.

\subsubsection{Estimation of variance of the EIF}

With the expression for the variance of the efficient influence function in hand
(Equation \ref{eq:eifvar2}), we can now form estimators which target this
parameter. $\bar{Q}_0$ and $g_0$ are not known in practice, though estimates
$\bar{Q}_n^{*}$ and $g_n$ will be readily available if estimating
$\mathbb{E}Y_{d}$ using a double robust estimator such as TMLE, thus providing
us with the observed outcome $S^{d}_{t}(\bar{Q}_n^{*},g_n)(\bar{L}(t))$.
Treating this variable as our new time point specific outcome, our goal is to
estimate the mean of this variable over the post-intervention distribution of
$\bar{L}^{d}(t)$. For notational convenience, let $Z^{d}(t) \equiv
S^{d}_t(\bar{Q}_0,g_0)(\bar{L}(t))$ represent the observable outcome and
$(L(0),A(0),\ldots,A(t-1),Z^{d}(t))$ represent the observed data structure.


One possible approach to estimating each of the components (Equation
\eqref{eq:varparm}) is to use a simple IPW estimator \citep{Horvitz1952} \[
\hat{\sigma}_{t,n,IPW}^{2,d}=\frac{1}{n}\sum_{i=1}^{n} \frac{\mathbb{I}(\bar{A}_
i(t^{-})=d(\bar{l}(t^{-})))}{g_{0:t^{-},n}(\bar{A}_i(t^{-}),\bar{L}_i(t^{-}))}
Z^{d}_n(t)\] where $Z^{d}_n(t)=S^{d}_t(\bar{Q}_n,g_n)(\bar{L}(t))$. However,
such an estimator would still be subject to underestimation of the variance by
ignoring the contribution of observations that selected a likely treatment
$\bar{A}_i$, even though their probability of following $d(\bar{l})$ is very
small. In other words, subjects $i$ with small probabilities of following
$d(\bar{l})$ would be unlikely to be observed with $\bar{A}_i=d(\bar{l})$
resulting in an indicator value of $0$ for the numerator and, consequently, a
contribution of $0$ to the IPW estimator. Therefore, we stress that it is
important to use a plug-in estimator such as TMLE \citep{vanderLaan2011a} to
estimate this parameter. A plug-in estimator will integrate over all
$\bar{l}(t)$ in the support of $P_{t,n}^{d}$ and thus contribute many large
values of $S^d_{t,n}(\bar{Q}_n^{*},g_n)$ when there are practical or theoretical
positivity assumption violations. In addition, the TMLE is a double robust
estimator so that it will yield a consistent estimator of this variance if $g_n$
is consistent for the true $g_0$.


Given $\bar{Q}_0,g_0$, we will now provide a succinct summary of the TMLE of
$\sigma^{2,d}_{0,t}=\mathbb{E}_{P_0^{d}}Z^{d}(t)$ that is based on iterative
sequential regression. Note that this iterative sequential regression approach
is similar to the one presented by van der Laan et al. \citep{vanderLaan2011a}
for the intervention specific mean outcome parameter. Denote the counterfactual
of $Z^{d}(t)$ under treatment $d'$ with $Z^{d,d'}(t)$, and let $P_0^{d'}$ be the
$G$-computation formula \citep{Robins1986} corresponding with this intervention
$\bar{A}(t^{-})=d'(\bar{l}(t^{-}))$. We wish to estimate
$\sigma^{2,d}_{0,t}=\mathbb{E}_{P_0^{d}}Z^{d,d}(t)$, which can be represented as
a series of iterated conditional expectations
\[ \sigma^{2,d}_{0,t}= \mathbb{E}[
\mathbb{E}[\cdots\mathbb{E}[\mathbb{E}[Z^{d}(t)|\bar{L}^{d}(t-1)]|\bar{L}^{d}(t-
2)]\cdots|\bar{L}^{d}(0)]]. \]
The EIF for this target parameter $\sigma^{2,d}_t$ is given by
\[
D_{\sigma^{2,d}_t}^{*}(P)(O)=\sum_{m=0}^{t} H^{d,t}_m(g)
(\bar{Q}^{d,\sigma_{t}^2}_{m^{+}}-\bar{Q}^{d,\sigma_{t}^2}_m),
\]
where we define
\begin{eqnarray*}
  \bar{Q}^{d,\sigma_{t}^2}_{t+1}&=&Z^{d}(t) \\
  H^{d,t}_m(g)&=&\frac{\mathbb{I}(\bar{A}(m^{-})=d(\bar{l}(m^{-})))}{g_{0:m^{-}}(\bar{A}(m^{-}),\bar{L}(m^{-}))}
  :m=1,2,\ldots,t \\
  H^{d,t}_0&=&1.
\end{eqnarray*}
Therefore, the EIF for $\sigma^2=\sum_t \sigma^{2,d}_t$ is simply
$D_{\sigma^2}^{*}=\sum_t D_{\sigma^{2,d}_t}^{*}$.

With the EIF established, the TMLE of $\sigma^{2,d}_t$ is now defined as
follows.
\begin{enumerate}
  \item Estimates $g_{0:m^{-},n}:m=1,2,\ldots,t$ are readily available
  if estimating $\mathbb{E}Y_{d}$ using a double robust estimator such as TMLE.
  \item Set $\bar{Q}^{d,\sigma^2_t}_{t,n}=Z^{d}_i(t)$. Determine the
  range $(a,b)$ for $Z^{d}_i(t)$, $i=1,\ldots,n$ and target this initial
  fit using a parametric submodel respecting this range
  $(a,b)$ by adding the clever covariate $H^{d,t}_t$ (on, say, the logistic
  scale), using the initial fit as off-set. The resulting updated fit is denoted
  with $\bar{Q}^{d,\sigma^2_t,*}_{t,n}$.
  \item Given $\bar{Q}^{d,\sigma^2_t,*}_{t,n}$, we can recursively for $m
  = t-1,t-2,\ldots,1$:
  \begin{enumerate}
    \item Regress the targeted fit $\bar{Q}^{d,\sigma^2_t,*}_{m^{+},n}$ onto
    $\bar{A}(m^{-})=d(\bar{l}(m^{-})),\bar{L}(m^{-})$, using logistic regression
    to respect the range $(a,b)$. Denote the fit
    $\bar{Q}^{d,\sigma^2_t}_{m,n}$.
    \item Target this initial fit respecting the range $(a,b)$ with clever
     covariate $\mathbb{I}(\bar{A}(m^{-})=d(\bar{l}(m^{-})))$ and observational
     weight $\frac{1}{g_{0:m^{-}}(\bar{A}(m^{-}),\bar{L}(m^{-}))}$ (on the
     logistic scale), and denote this targeted fit of
     $\bar{Q}^{d,\sigma^2_t}_{m}$ with $\bar{Q}^{d,\sigma^2_t,*}_{m,n}$.
  \end{enumerate}
  \item At $m=1$, we have the estimate $\bar{Q}^{d,\sigma^2_t,*}_{1,n}$,
  which now is a function of only $L(0)$. Finally, we take the average of
  $\bar{Q}^{d,\sigma^2_t,*}_{1,n}$ w.r.t. the empirical distribution of $L_i(0)$.
  The resulting $\hat{\sigma}^{2,d}_{t,n,TMLE}=\bar{Q}^{d,\sigma^2_t,*}_{0,n}$ is
  the desired TMLE of $\sigma^{2,d}_t$.
\end{enumerate}

\subsubsection*{Estimation of variance of the EIF}

\subsubsection{Application to single time-point treatment setting}

For the sake of illustration, let us consider the method presented above for
estimation of the variance of the EIF for the case that $O=(L(0),A(0),Y=L(1))$
and the target parameter is $\mathbb{E}Y^a$ for a static point treatment $a$.

In this case, the variance of the efficient influence curve is represented as
\begin{equation}
\label{eq:pteicvar}
\begin{split}
  \sigma^2_0 &= \mathbb{E}_0[D^{*}(P_0)(O)]^2 \\
    &=
    \mathbb{E}_0\left[\frac{\mathbb{I}(A=a)}{g_{0}(a
    \mid L(0))}(Y-\bar{Q}^a_0(L(0))) + \bar{Q}^a_0(L(0)) - \mathbb{E}Y^a\right]^2
    \\
    &=\mathbb{E}_0\left[\frac{\mathbb{I}(A=a)}{g_{0}(a \mid
    L(0))}(Y-\bar{Q}^a_0(L(0)))\right]^2 + \mathbb{E}_0[\bar{Q}^a_0(L(0)) - \mathbb{E}Y^a]^2 \\
    &=\mathbb{E}_{P_0^a}\left[\frac{(Y^a-\bar{Q}^a_0(L(0)))^2}{g_{0}(a \mid L(0))}\right] +
    \mathbb{E}_{0}[\bar{Q}^a_0(L(0))-\mathbb{E}Y^a]^2.
\end{split}
\end{equation}

If using a double robust estimator for the estimation of $\mathbb{E}Y^a$ such as
TMLE, we are provided with estimators $g_n$ and $\bar{Q}_n^*$ of $g_0(A \mid
L(0))$ and $\bar{Q}_0^a(L(0))=\mathbb{E}[Y^a \mid L(0)]=\mathbb{E}_0[Y \mid
A=a,L(0)]$ respectively. The second term in the final expression of Equation
\eqref{eq:pteicvar} is easily estimated with the empirical distribution. Given
$g_0$ and $\bar{Q}_0$, the first term can be represented as the mean of a
counterfactual $S^a(L^{a}(0)) \equiv (Y^a-\bar{Q}_0^a(L(0)))^2/ g_0(a \mid
L(0))$ which needs to be estimated based on $(L(0),A,S^a(L(0),Y))$, where
$S^a(L(0),Y)=(Y-\bar{Q}_0(a,L(0)))^2/g_0(a \mid L(0))$ represents the observed
outcome. For example, we can use a TMLE estimator
$\mathbb{E}_{n}^*S^a(L(0),Y^a)$ of
$\mathbb{E}_0S^a(L(0),Y^a)=\mathbb{E}_{L(0),0}[\mathbb{E}_0[S^a \mid
A=a,L(0)]]$. The TMLE estimate $\mathbb{E}_n^*[S^a\mid A=a,L(0)]$ of
$\mathbb{E}_0[S^a\mid A=a,L(0)]$ is defined by determining the range $(a,b)$ of
$S^a(L_i(0),Y_i)$, obtaining an initial regression fit of $\mathbb{E}_0[S^a\mid
L(0),A]$ that respects this range, representing it as a logistic regression fit
bounded by $(a,b)$, and updating the latter by fitting a univariate logistic
regression with clever covariate $\mathbb{I}(A=a)$ and observational weight
$1/g_0(a \mid L(0))$, using the initial fit as an off-set. Regarding the initial
fit $\mathbb{E}_n[S^a \mid A=a,L(0)]$, recall from above that $S^a$ is a
function of $L(0)$ which results in the initial fit being exactly
$(Y-\bar{Q}_0(a,L(0)))^2/g_0(a \mid L(0))$ such that regression is unneeded.
Following the update step, the TMLE of $\mathbb{E}_0S^a(L(0),Y^a)$ is now given
by $\frac{1}{n}\sum_{i=1}^n \mathbb{E}_n^*[S^a\mid L_i(0),A=a]$, so that \[
\hat{\sigma}^{2,*}_n=\frac{1}{n}\sum_{i=1}^n \mathbb{E}_n^*[S^a(\bar{Q}_n^*,g_n)
\mid A=a,L_i(0)] + \frac{1}{n}\sum_{i=1}^n
(\bar{Q}_n^*(L_i(0),a)-\hat{\psi}_n^*)^2 \] where $\hat{\psi}_n^*$ is the
targeted estimate of $\mathbb{E}Y^a$.



\subsection{Advantages of this plug-in estimator of the asymptotic variance of
the EIF}

Since $\sigma^2_0/n$ equals the asymptotic variance of an asymptotically
efficient estimator, it provides a good measure of the amount of information in
the data for the target parameter of interest. Therefore, it is sensible to view
$\sigma^2_0/n$ as a measure of sparsity for the target parameter of interest. If
$g_n$ is a good estimator of $g_0$, then our proposed plug-in estimator
$\hat{\sigma}^2_n$ is much less subject to under-estimation due to sparsity than
currently available estimators such as the sample variance of the estimated
influence function, and the bootstrap-based estimate of the variance of an
efficient estimator. Indeed, the non-parametric bootstrap generally is not
valid, except when using a parametric model to estimate $g_0$ and $\bar{Q}_0$
which will never capture a true model in practice. This plug-in estimate
$\hat{\sigma}^2_n$ represents a variance of the estimate of the EIF which
involves the integration of rare combinations of treatment and covariates that
are unlikely to occur in the actual sample.

In particular, if there are theoretical violations of the positivity assumption,
then this true variance $\sigma^2$ equals infinity, and, if $g_n$
approximates $g_0$ well, then also the estimate $\hat{\sigma}^2_n$ will generate
very large values, demonstrating the lack of identifability and thereby raising
a red flag for finite sample sparsity bias in the estimators (beyond the large
confidence intervals generated by $\hat{\sigma}^2_n$). If there are
serious practical violations of the positivity assumption, then the estimate of
this variance should be imprecise, since it is itself a highly variable
estimator of a weakly identifiable parameter.

\section{Variance estimation for substitution based estimators}
\label{sec:bootstrap}

The plug-in estimator of the asymptotic variance of the EIF presented above is
superior to the more common approach of taking the empirical EIF variance over
the sample (i.e., $\text{var}[D^{*}(P_n)(O)]$), in that there is a much stronger
contribution of combinations of treatment and covariates that are unlikely to
occur in the actual sample. In finite samples, however, the use of substitution
based estimators such as TMLE (which are guaranteed to solve the EIF within a
bounded range) are often observed to have smaller variance than their asymptotic
variance. This is due to the mere fact that they are guaranteed to respect the
global constraints of the statistical model and target parameter mapping. That
is, as opposed to estimating equations that tend to result in estimates outside
parameter boundaries as the EIF variance increases, the use of substitution
estimators in finite samples will retain an estimator variance that is smaller
than the EIF variance divided by the sample size, $n$. Thus, using the newly
presented robust EIF variance method can result in over-estimation of the
estimator variance for these types of estimators. This therefore motivates us to
develop an estimator of variance that is less conservative, i.e. more aligned to
the true variance of substitution based estimators (such as TMLE).

One alternative approach at estimating the variance for substitution based
estimators is to conduct a non-parametric bootstrap. The $n$ observations are
sampled with replacement and used to form an estimate of the parameter over $B$
iterations. However, as stated above, the non-parametric bootstrap is generally
invalid and not theoretically supported. Additionally, this is a very computer
intensive method that usually requires estimating the full likelihood (i.e.,
$P_0^{d}$) of the longitudinal data structure within each sampled iteration and
is therefore normally infeasible in practice unless conducted within an a priori
selected smaller parametric statistical model such as logistic regression.

In this section we present an alternative bootstrap based approach that, unlike
the standard non-parametric bootstrap, is both computationally feasible and
theoretically valid. That is, this bootstrap approach allows us to estimate the
variance of the estimator while avoiding re-estimation of $g_0$ and $\bar{Q}_0$.
To facilitate this, we propose a modification of the usual TMLE such that the
targeting step is separated from the initial estimation of $\bar{Q}_0$. Recall
that the typical TMLE, as implemented, pivots between the targeting step and the
initial estimator for the next regression (preventing us from separating the
intitial fit from the targeting step). We propose a minor modification of the
TMLE that separates these steps, first estimating all of the initial regressions
and subsequently targeting the fits in a separate step. This modified TMLE can
then be boostrapped via only the targeting step. Note that, because the modified
TMLE has the same asymptotic behavior as the original TMLE, the bootstrap is
theoretically supported and will lead to valid inference.

\subsection{Modified TMLE for \texorpdfstring{$\mathbb{E}Y^{d}$}{Lg}}

To reduce the computational burden that bootstrapping requires, we first present
the modified TMLE approach at estimating the parameter $\mathbb{E}Y^{d}$. This
parameter can be estimated by the following steps:

\begin{enumerate}
  \item Estimate $g_{0:t^{-}}(\bar{A},\bar{L}):t=1,2,\ldots, K+1$ and denote the
  fits $g_{0:t^{-},n}$.
  \item Determine the range $(a,b)$ for $\mathbb{E}Y^{d}$. Recursively
  for $t=K+1,K,\ldots,1$, estimate the conditional expectation
  $\bar{Q}_{t}^{d} = \mathbb{E}[\bar{Q}_{t^{+}}^ {d}|\bar{L}(t^{-}),
  \bar{A}(t^{-})=d(\bar{l}(t^{-}))]$ respecting this range. Denote the fits
  $\bar{Q}_{t,n}^{d}$. \textit{We stress that this step is cruicially different than the
    typical TMLE, in that all of the initial regression fits are done simultaneously.}
  \item For time $t=K+1$, target the initial fit $\bar{Q}_{K+1,n}^{d}$ by
  using a parametric submodel respecting the range $(a,b)$ by adding the
  covariate $\mathbb{I}(\bar{A}(K) = d(\bar{l}(K)))$ and observational weight
  $1/g_{0:K,n}$ (on the logistic scale), using the initial fits as off-set, and
  setting $Y$ as the dependent variable. Denote this updated fit as
  $\bar{Q}_{K+1,n}^{d,*}$.
  \item Given $\bar{Q}_{K+1,n}^{d,*}$, we can recursively for
  $t=K,K,K-1,\ldots,1$ target the initial fits $\bar{Q}_{t,n}^{d}$ by
  using parametric submodels respecting the range $(a,b)$, adding the covariates
  $\mathbb{I}(\bar{A}(t^{-}) = d(\bar{l}(t^{-})))$ and observational weight
  $1/g_{0:t^{-},n}$ (on the logistic scale), using the initial fits as off-set,
  and setting $\bar{Q}_{t^{+},n}^{d,*}$ as the dependent variable.
  Denote the updated fits as $\bar{Q}_{t,n}^{d,*}$.
  \item At $t=1$, we have the estimate $\bar{Q}_{1,n}^{d,*}$,
  which now is a function of only $L(0)$. Taking the average of
  $\bar{Q}_{1,n}^{d,*}$ w.r.t. the empirical distribution of $L_i(0)$
  gives us the desired TMLE estimate of $\mathbb{E}Y^{d}$.
\end{enumerate}

This estimator also solves the EIF and is therefore also asymptotically linear
and efficient. We note that the analysis of this TMLE is identical to the
typical TMLE presented by van der Laan et al. \citep{vanderLaan2011a}, with the only difference
being the initial estimator fits. Here the initial estimators are the original
ones, whereas the previous TMLE is implemented with initial estimators using the
targeted fits for the outcome.

We emphasize that this estimator is proposed for the sake of the bootstrap
method for variance estimation. It is recursive, in that each fit
$\bar{Q}_{t,n}^{d}$ is dependent upon the fit at $t^{+}$. As opposed to the
TMLE, the recursive nature of this TMLE is self contained within each step. In
other words, each estimation step in this TMLE can be performed independently of
the other steps. This allows the analyst to form all of the initial fits $P_{n}$
prior to performing any of the targeted updates.

\subsection{Bootstrapping the modified TMLE}

The new TMLE approach presented above can be bootstrapped in a fully
non-parametric manner, such that observations are drawn with replacement prior
to fitting the full likelihood $P_0^{d}$ and used to form an estimates of the
parameter, leading to an estimate of estimator variance. Our recommendation
is to only bootstrap the targeting step. More specifically, once the fits
$g_{0:t^{-},n}$ and $\bar{Q}_{t,n}^{d}$ are formed for $t=1,2,\ldots,K+1$, steps
3-5 above are carried out in the bootstrap such that for $b=1,2,\ldots,B$ we
have \[ Q_{n,b}^{*}=Q_n(\epsilon_b) \] for a user selected submodel
$P(\epsilon)$. The estimator variance is then estimated by taking the variance
over the bootstrapped estimates, i.e., $\text{var}(\hat{\Psi(Q_n)}) =
\text{var}[\Psi(Q_{n,b}^{*})]$.

We emphasize that this TMLE is provided such that we do not need to re-estimate
$\bar{Q}_n, g_n$, i.e. it is not a function of the data. If $g_n \rightarrow
g_0$ and $\bar{Q}_n \rightarrow Q$, then this TMLE is asymptotically linear with
influence function $D*(Q, g_0)$. This is conservative relative to the variance of
the actual TMLE that is estimated with $g_n$ fitted on the data, when $g_n$ is
consistent.

\section{Simulations}
\label{sec:sim}

Simulation studies presented in this section illustrate the performance of the
two estimators of variance for the estimation of the effect of treatment in both
a point treatment setting, and in a longitudinal observational study setting
with three time points (i.e. $K+1=3$) with time-dependent confounding. To
analyze the performance, we first compare the variance estimation approaches
covered above in estimating the estimator variance. Both the AIPW and TMLE
estimators are considered in order to demonstrate the difference in estimating
equations and substitution based estimators, respectively. The mean of the
variance estimates are compared to the Monte-Carlo variance of each estimator.
Additionally, we present the empirical coverage, Type I, and Type II errors
resulting from each variance estimation approach. The Monte-Carlo variance of
each variance approach is also reported. All analyses were conducted on R
version 3.1.1 \citep{Team2014}.

\subsection{Data generating distribution \texorpdfstring{$P_0$}{Lg}}

\subsubsection{Point treatment setting}
\label{sim:pt}

Consider a point treatment setting, such as patient enrollment into a care
program, in which the treatment $A(0)$ is only assigned at a single time point.
We are interested in determining whether the treatment of interest has a
significant effect on the outcome on an additive scale. Our target parameter is
therefore the difference of the mean outcomes under treatment and control, i.e.,
$\psi_{0,1} \equiv \mathbb{E}\,Y_1 - \mathbb{E}\,Y_0$. Under this setting, the
simulated data were generated as follows:

\begin{eqnarray*}
  W_1, W_3 &\sim& N(0,1) \mbox{, bounded at [-2,2]} \\
  W_2 &\sim& Ber(\text{logit}^{-1}(-1)) \\
  L_1(0) &\sim& N(0.1 + 0.4W_1, 0.5^2) \\
  L_2(0) &\sim& N(-0.55 + 0.5W_1 + 0.75W_2, 0.5^2) \\
  \bar{g}_{0,0}(Pa(A(0))) &=& \text{logit}^{-1}(\beta_{p} - (\beta_{p}+2.5)W_1 +
  1.75W_2 \\
    &\hspace{.1in}& + (\beta_{p}+3.2)L_1(0) - 1.8L_2(0) + 0.8L_1(0)L_2(0))) \\
  \bar{Q}_{0,1}(Pa(Y)) &=& \text{logit}^{-1}(-0.5 + 1.2W_1 - 2.4W_2 - 1.8L_1(0) -
  1.6L_2(0) \\
    &\hspace{.1in}& + L_1(0)L_2(0) - \beta_{\psi_0}A(0))
\end{eqnarray*}

with a positivity associated parameter $\beta_{p}$ ranging from $-2$
(minor positivity violations) to $1$ (strong practical positivity
violations) and the treatment effect associated parameter $\beta_{\psi_0}$
ranging from $0$ (no treatment effect) to $1$ (strong treatment effect). Here,
$L_1(0)$ and $L_2(0)$ are not time-dependent confounders and are therefore
considered baseline covariates along with $(W_1,W_2)$, which affect both the
treatment and the outcome.

\subsubsection{Longitudinal treatment setting}

For the longitudinal setting, we considered a treatment $A(t)$ which was allowed
to vary over time as a counting process. That is, if $A(t)=1$ then we have that
$\underline{A}(t)=1$ where $\underline{X}(t)=(X(t), X(t+1), \ldots, X(K))$.
Similar to the point treatment setting, we are interested in whether the
treatment of interest has a significant effect on the outcome at the final time
point $t^{*}=3$ on an additive scale. Thus, our target parameter is the
difference of the mean outcomes under treatment and control at this final time
point, i.e., $\psi_{0,3}=\mathbb{E}\,Y_1(t^{*}) - \mathbb{E}\,Y_0(t^{*})$ where
$Y(t^{*})=L_3(3)$. Under this setting, data for the first time point was
generated in the same manner as the point treatment setting in Section
\ref{sim:pt} above. For the remaining two time points, the data were generated
conditional on survival (i.e.
$L_3(t^{-})=0$) as follows:

\begin{eqnarray*}
  L_1(t) &\sim& N(0.1 + 0.4W_1, 0.5^2 + 0.6L_1(t^{-}) - 0.7L_2(t^{-}) +
  0.45\beta_{\psi_0}A(t^{-})) \\
  L_2(t) &\sim& N(-0.55 + 0.5W_1 + 0.75W_2 + 0.1L_1(t^{-}) + 0.3L_2(t^{-}) \\
  &\hspace{.1in}& + 0.75\beta_{\psi_0}A(t^{-}), 0.5^2) \\
  \bar{g}_{0,t}(Pa(A(t))) &=& \text{logit}^{-1}(\beta_{p} - (\beta_{p}+2.5)W_1 +
  1.75W_2 \\
    &\hspace{.1in}& + (\beta_{p}+3.2)L_1(t) - 1.8L_2(t) + 0.8L_1(t)L_2(t))) \\
  \bar{Q}_{0,t}(Pa(L_3(t))) &=& \text{logit}^{-1}(-0.5 + 1.2W_1 - 2.4W_2 -
  1.8L_1(t^{-}) - 1.6L_2(t^{-}) \\
    &\hspace{.1in}& + L_1(t^{-})L_2(t^{-}) - \beta_{\psi_0}A(t^{-}))
\end{eqnarray*}

\begin{figure}[!h]
\centering
\includegraphics[width=2in]{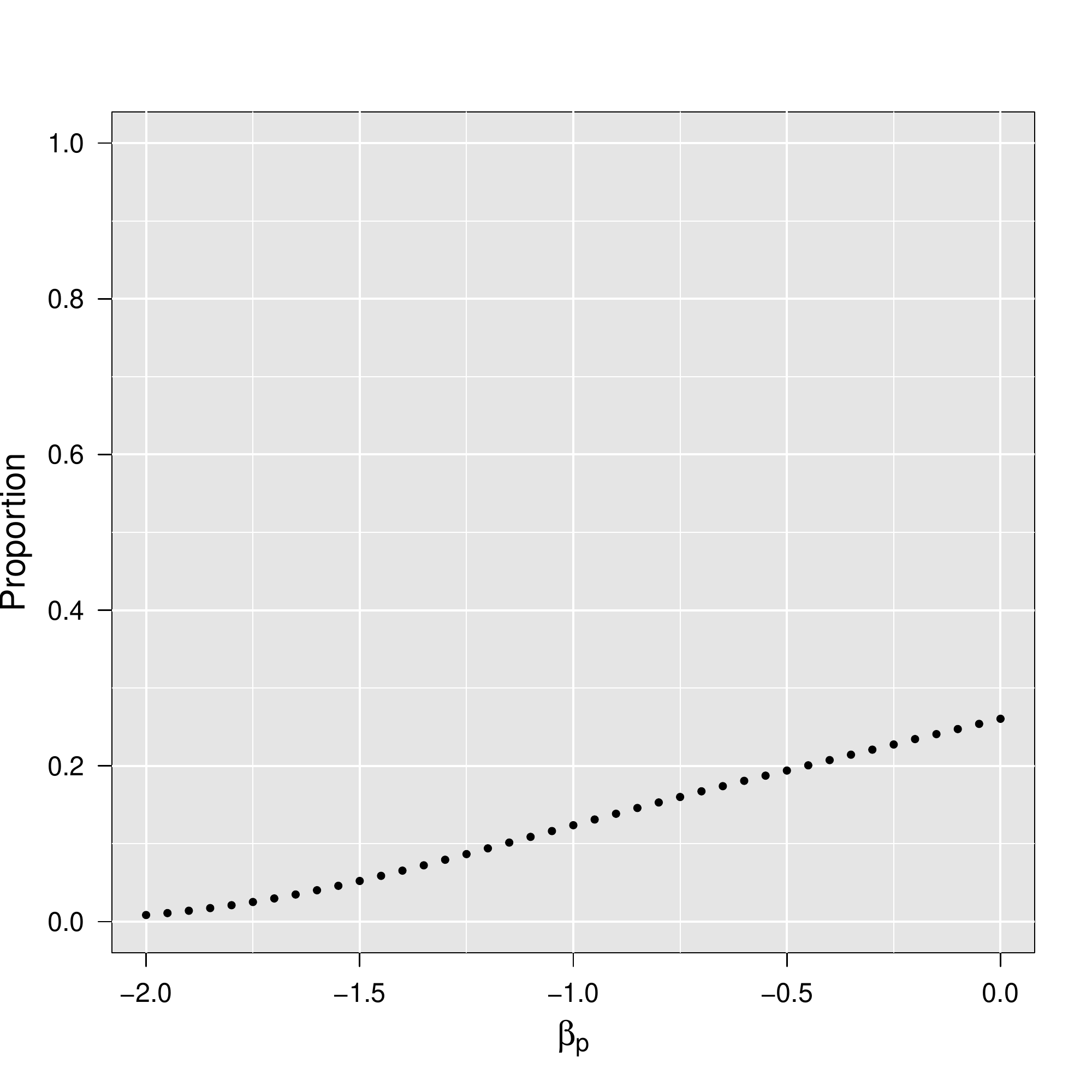}
\caption{Proportion of observations with $g_{0:2}$ truncated at each $\beta_{p}$.}
\label{fig:truncation}
\end{figure}

Similar to the point treatment setting, the treatment effect associated
parameter $\beta_{\psi_0}$ also ranged from $0$ to $1$. We note, however, that
the positivity issues faced in this scenario will be even more severe due to the
higher number of combinations of treatment over time, which result in smaller
probabilities. We therefore considered only $\beta_{p}$ values from $-2$ to $0$
and imposed a truncation level of $0.001$ to the estimates of $g_{0:t}$. Figure
\ref{fig:truncation} shows the proportion of observations with truncated
$g_{0:2}$ as a function of $\beta_{p}$ at a null effect, i.e. $\beta_{\psi_0} = 0$.

Under these settings, the true parameter values $\psi_0$ were achieved by
generating $8 \times 10^7$ observations under the counterfactual distribution
for each $\beta_{\psi_0}$ considered. Simulation results were obtained for $500$
simulations of size $n = 500$. Within each simulation, the bootstrap estimates
of variance were formed from $B=1000$ iterations.

\subsection{Submodels used}

Any submodel and loss function for which its loss-function specific score
\[\frac{\partial}{\partial\epsilon} L(P(\epsilon))\Big|_{\epsilon=0}\] spans
$D^{*}(P_0)$ can be chosen in TMLE for both estimation of the mean outcome
$\mathbb{E}Y_a$ and the variance of the EIF $\sigma^2$. As these submodels solve
the equation corresponding to the EIF, they will all be asymptotically
equivalent and thus, all be asymptotically efficient. That is, no difference will
be seen between the use of various submodels as the sample size grows to
infinity. In the TMLE presented by van der Laan et al. \citep{vanderLaan2011a}, these submodels are
used in the targeting step for each $\bar{Q}_t$ using a loss $L(\bar{Q}_t)$
that is indexed by $\bar{Q}_{t+1}$. Specifically, for the targeting step we need
a loss and submodel with clever covariate such that the score given solves a
desired component of the efficient influence function $D^{*}_{t}(P_0)$.

The use of various submodels in finite samples can have varying performance. For
example, under increasing levels of positivity violations the use of linear
submodels which use $H_t(g)$ as a covariate can have higher variance due to
observations with low probabilities of treatment acting as outliers which result
in highly influential points for the estimation of the submodel parameter
$\epsilon$.

Recall that the catalyst for this work was the anti-conservative estimates of
estimator variance resulting from the use of the empirical EIF variance. We
therefore wish to establish a robust estimator of the variance of estimators
which solve the EIF, particularly under violations or near violations of
positivity. In other words, we desire a variance estimator which will
asymptotically converge to the true variance of the estimator, but also
simultaneously act on the conservative side in finite samples. Keeping this in
mind, we used two submodel and loss function combinations for our simulations.
For the estimation of the target parameter and the robust estimator of the EIF
variance, we used submodels which define $H_t(g)$ and $H_m^{d,t}(g)$ to be
observational weights such that \[ \text{logit}\bar{Q}(\epsilon) =
\text{logit}\bar{Q} + \epsilon, \] acknowledging our slight abuse of notation.
Alternatively, in our bootstrap approach at estimating the TMLE variance, we
define a clever covariate using $H_t(g)$ such that \[
\text{logit}\bar{Q}(\epsilon) = \text{logit}\bar{Q} + \epsilon H_t(g).\] Both
submodels use, as loss function, the negative log-likelihood loss. As stated
previously, both of the submodels presented solve the equation corresponding to
the EIF and are therefore asymptotically equivalent.

\subsection{Simulation results}

\subsubsection{Point treatment results}

\begin{figure}[!h]
\centering
\includegraphics[width=4in]{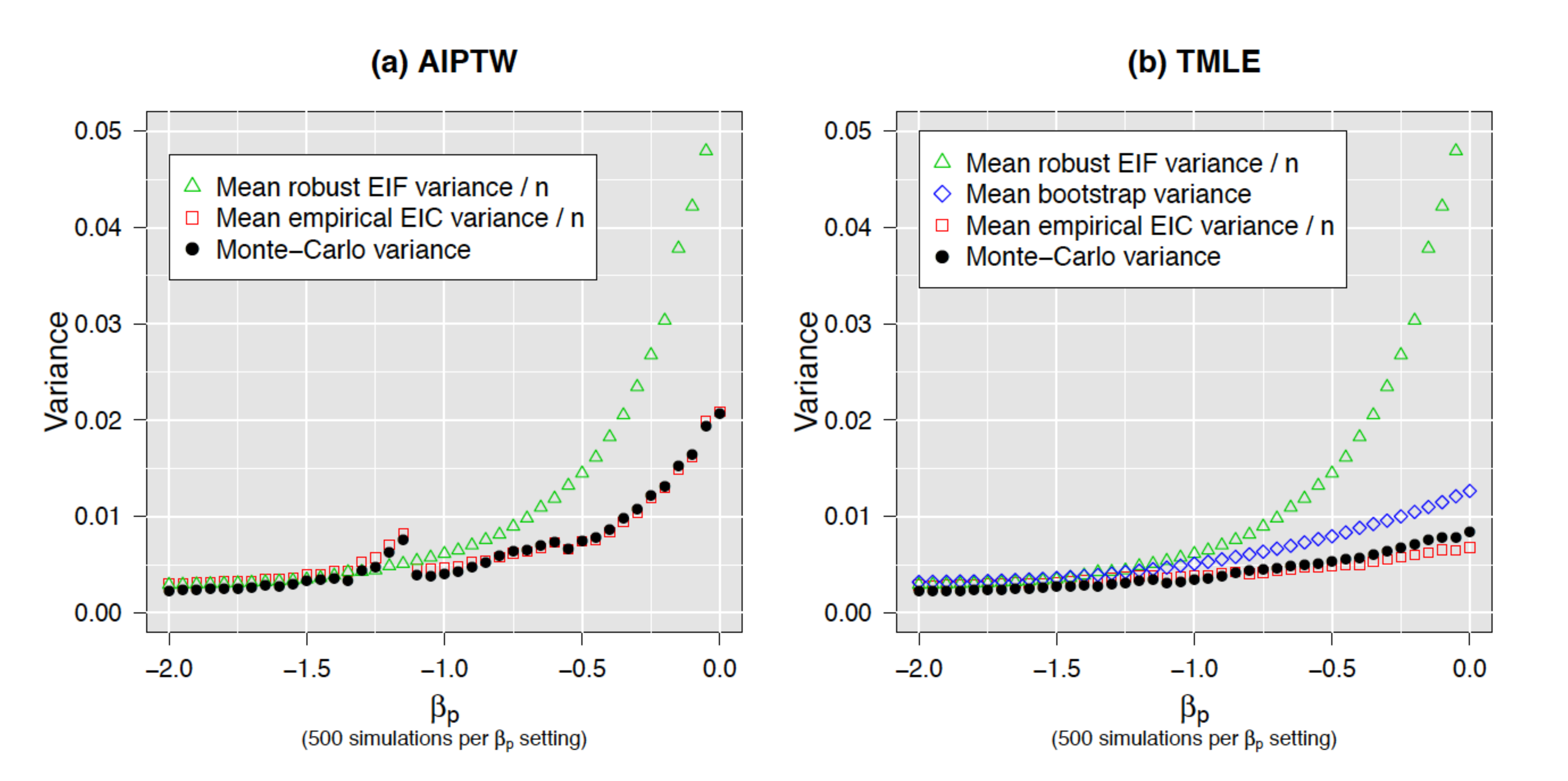}
\caption{Mean of variance estimates for each estimator under no treatment
effect ($\beta_{\psi_0}=0$) at each positivity ($\beta_{p}$) value under the
point treatment setting, overlaid with the estimator's Monte-Carlo variance.
Robust variance estimates are identical for the two estimators.}
\label{fig:compare}
\end{figure}

Figure \ref{fig:compare} shows the Monte-Carlo variance under no treatment
effect ($\beta_{\psi_0}=0$) for both the AIPW and TMLE estimators, along with
the mean of the variance estimates from each estimation approach. To keep the
differences in perspective, we plotted results only for the positivity
associated parameter $\beta_{p} \le 0$. At the lower end of $\beta_{p}$ where
positivity violations are minor, the observed estimator variance is low for both
the AIPW and TMLE estimators, with the TMLE approach showing lower variance
between the two despite solving the same estimating equation corresponding to
the EIF. For example, at $\beta_{p}=-2$ the Monte-Carlo variance was $0.0024$
and $0.0022$ for the AIPW and TMLE estimators, respectively. As $\beta_{p}$
increased, introducing higher levels of positivity violations, the estimator
variance increased for both estimators. Additionally, this occurred at a much
higher rate for the AIPW estimator than for TMLE, resulting in an increase in
the magnitude of difference between the two estimators. For example, at
$\beta_{p}=0$ the simulations resulted in a Monte-Carlo variance of $0.0207$ and
$0.0085$ for the AIPW and TMLE estimators, respectively.

For the AIPW estimator, the empirical EIF based approach of estimating variance
performed resulted in estimates similar to the Monte Carlo estimates. For
example, at $\beta_{p}=0$ the mean of the EIF approach was $0.0208$. A slight
but consistent underestimation of the variance was observed at higher levels of
practical positivity violations. The robust approach of estimating variance
appeared to result in conservative estimates of variance.

In the TMLE estimator, all three approaches to variance estimation performed
similarly at low values of $\beta_{p}$. For example, at $\beta_{p}=-2$, the mean
of the estimates was $0.0029$, $0.0029$, and $0.0032$ for the empirical EIF,
robust, and bootstrapped based approaches, respectively, compared to the
estimator's Monte-Carlo variance of $0.0022$. As $\beta_{p}$ increased, the
empirical EIF approach tended to result in anti-conservative estimates of
variance, while the bootstrap approach resulted in slightly conservative
estimates. The robust EIF approach tended to overestimate the TMLE estimator
variance.

\begin{figure}[!h]
\centering
\includegraphics[width=4in]{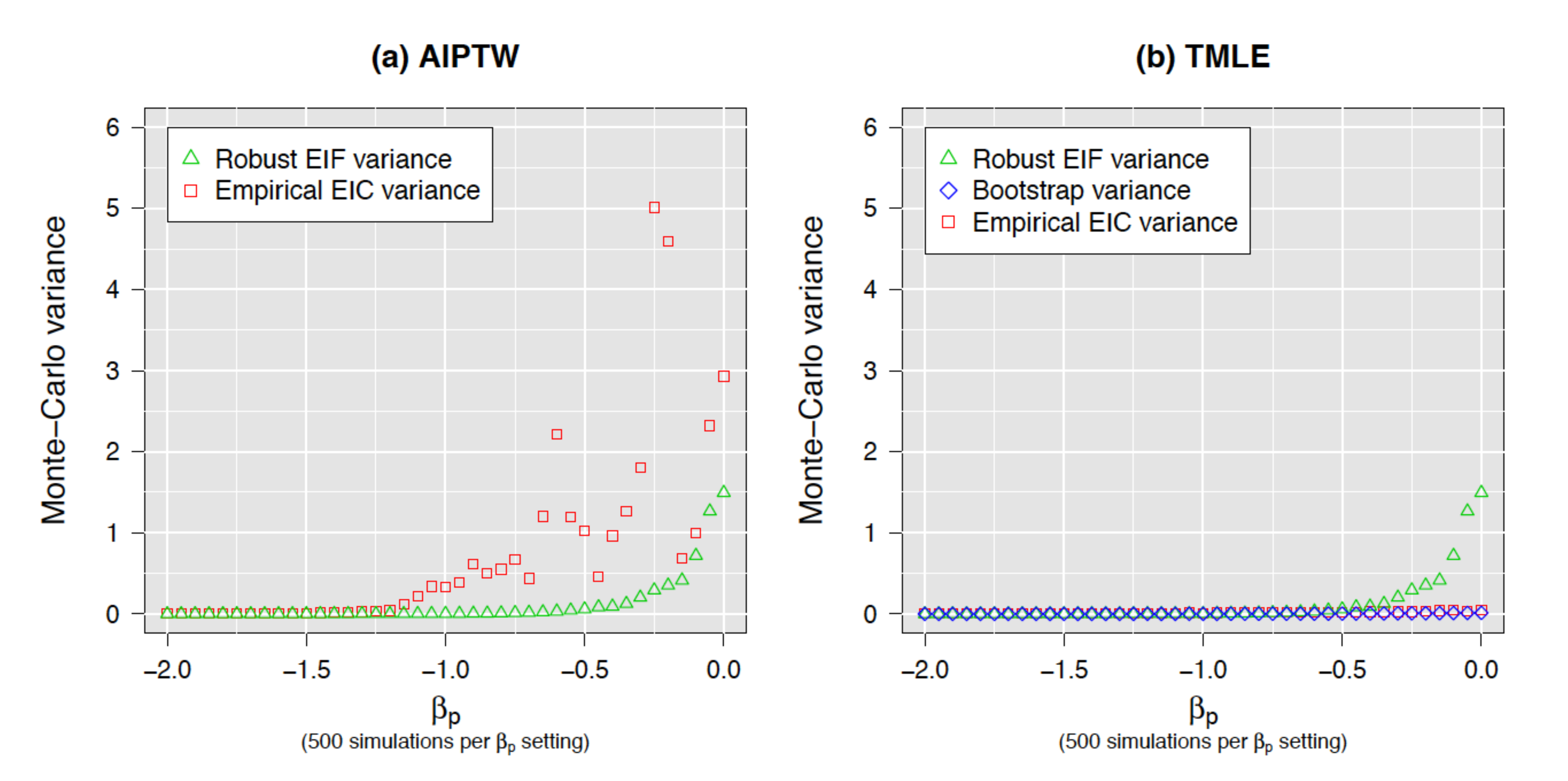}
\caption{Monte-Carlo variance of variance estimators for each mean
outcome estimator under no treatment effect ($\beta_{\psi_0}=0$) at each
positivity ($\beta_{p}$) value under the point treatment setting. Robust
variance estimates are identical for the two mean outcome estimators.}
\label{fig:compareVar}
\end{figure}

Figure \ref{fig:compareVar} shows the Monte-Carlo variance for each approach
taken at estimating the variance. Lower values in this figure can be interpreted as
coming from a variance estimator with more precision. In the AIPW estimator, the
empirical EIF approach has noticeably higher variance than the robust approach,
with a variance of $2.93$ at $\beta_{p}=0$ compared with $1.50$ for the robust
approach. This implies that the empirical EIF approach to estimating the AIPW
estimator variance is less reliable than the robust EIF approach. For the TMLE
estimator (Figure \ref{fig:compareVar}b), the empirical EIF approach to estimating
variance showed much lower Monte-Carlo variance. The bootstrap approach also
resulted in very low variance, implying a high reliability of this approach at
estimating the variance.

\begin{figure}[!ht]
\centering
\includegraphics[width=5in]{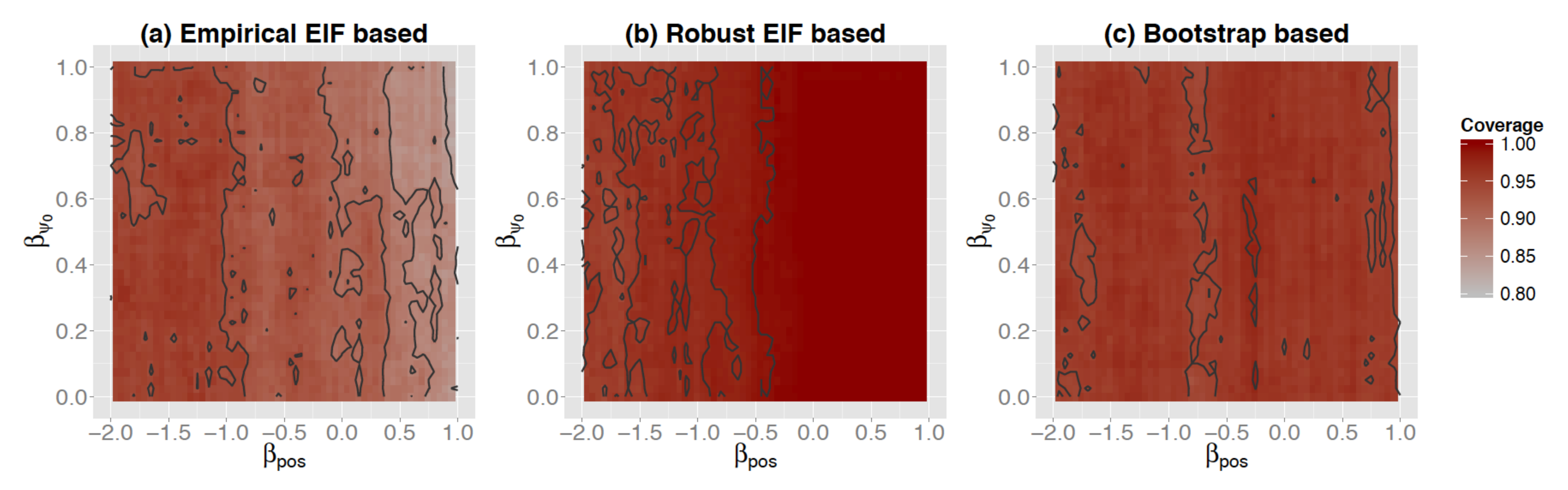}
\caption{Simulated coverage for each variance estimation approach for the
TMLE estimator under various treatment ($\beta_{\psi_0}$) and positivity
($\beta_{p}$) values under the point treatment setting.}
\label{fig:t1cover}
\end{figure}

\begin{figure}[!h]
\centering
\includegraphics[width=5in]{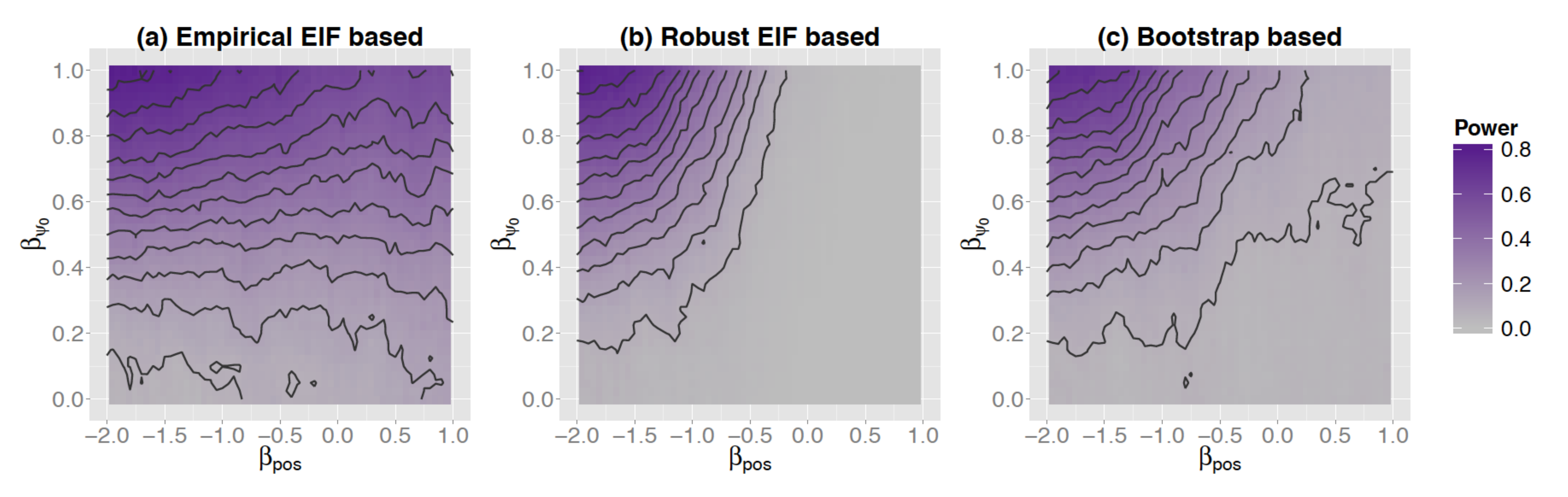}
\caption{Simulated power for each variance estimation approach for the TMLE
estimator under various treatment ($\beta_{\psi_0}$) and positivity
($\beta_{p}$) values under the point treatment setting.}
\label{fig:t1power}
\end{figure}

We evaluated $95\%$ confidence interval coverage for the TMLE estimator of
$\mathbb{E}Y_{d}$ under the three approaches to variance estimation. Due to
the lower variance seen in Figures \ref{fig:compare} and \ref{fig:compareVar}, we
focused only on the TMLE estimator here. Figure \ref{fig:t1cover} shows a heat
map overlaid with a contour plot of the resulting coverage estimates (i.e. the
observed proportion of times the true parameters were captured by the confidence
intervals) over the different combinations of $\beta_{\psi_0}$ and $\beta_{p}$.
Additionally, we estimated the power to reject the null hypothesis (at a level
of $0.05$) corresponding to each variance estimation approach under the range of
treatment effect sizes and degrees of positivity violation considered above.
Figure \ref{fig:t1power} shows a heat map overlaid with a contour plot of the
resulting power estimates. Results at $\beta_{\psi_0}=0$ can be interpreted as
Type I errors, as they inform us of the times that the null hypothesis of no
treatment effect is incorrectly rejected.

At low instances of positivity issues, coverage appears valid for all three
variance estimation approaches with the proportion of time the true parameter
was captured consistently at $0.95$ or larger (Figure \ref{fig:t1cover}).
Where positivity issues were low ($\beta_{\psi_0}< -0.5$), the empirical EIF
approach maintained nominal to conservative coverage. Where severe positivity
violations were present, coverage dropped substantially below $0.95$. For
example, at $\beta_{p}=1$ coverage for this approach varied from $0.41$ to
$0.85$. In contrast, the robust EIF approach consistently resulted in coverage
at around $0.95-0.96$ at low values of $\beta_{p}$ and \textit{increased} with
$\beta_{p}$, consistent with prior results showing overestimation of the
variance under increasing positivity by this approach. For example, at
$\beta_{p}=-.7$, coverage remained at $0.98$ at all values of
$\beta_{\psi_0}$. At $\beta_{p}\ge -0.1$, the observed coverage was almost
always greater than or equal to $0.99$ at all values of $\beta_{\psi_0}$.
The bootstrap based coverage shown in Figure \ref{fig:t1cover}c varied the
least, with coverage consistently between $0.95-0.97$ irrespective of the
treatment effect ($\beta_{\psi_0}$) and positivity severity ($\beta_{p}$)
considered.

Regarding the observed power (Figure \ref{fig:t1power}), the empirical-EIF based
variance approach resulted in the highest power among all three variance
estimation approaches when an effect was present. For example, at
$\beta_{\psi_0}=1$ and $\beta_{p}=-1$, the observed power was $0.71$, $0.51$,
and $0.51$ for the empirical-EIF, robust-EIF, and bootstrap approaches
respectively. While this result implies a more efficient approach, it expectedly
came at a cost of higher Type I error which became uncontrolled with an increase
in $\beta_{p}$. For example, at $\beta_{p}=-2$ an observed $4.2\%$ of the
simulations incorrectly rejected the null hypothesis. This proportion increased
to as high as $15\%$ at $\beta_{p}=1$. Alternatively, the robust EIF estimation
approach resulted in low Type I errors (i.e. between $0-5.8\%$) with none of the
simulations incorrectly rejecting the null beyond $\beta_{p}=-0.1$. The
bootstrap approach resulted in an intermediate performance, with higher power
than the robust EIF approach when an effect was present while simultaneously
retaining appropriate control of the Type I error at all levels of $\beta_{p}$
when no effect was present. For example, at $\beta_{p}=1$ only $4.8\%$ of the
simulations incorrectly rejected the null hypothesis.

\subsubsection{Longitudinal treatment results}

\begin{figure}[!h]
\centering
\includegraphics[width=4in]{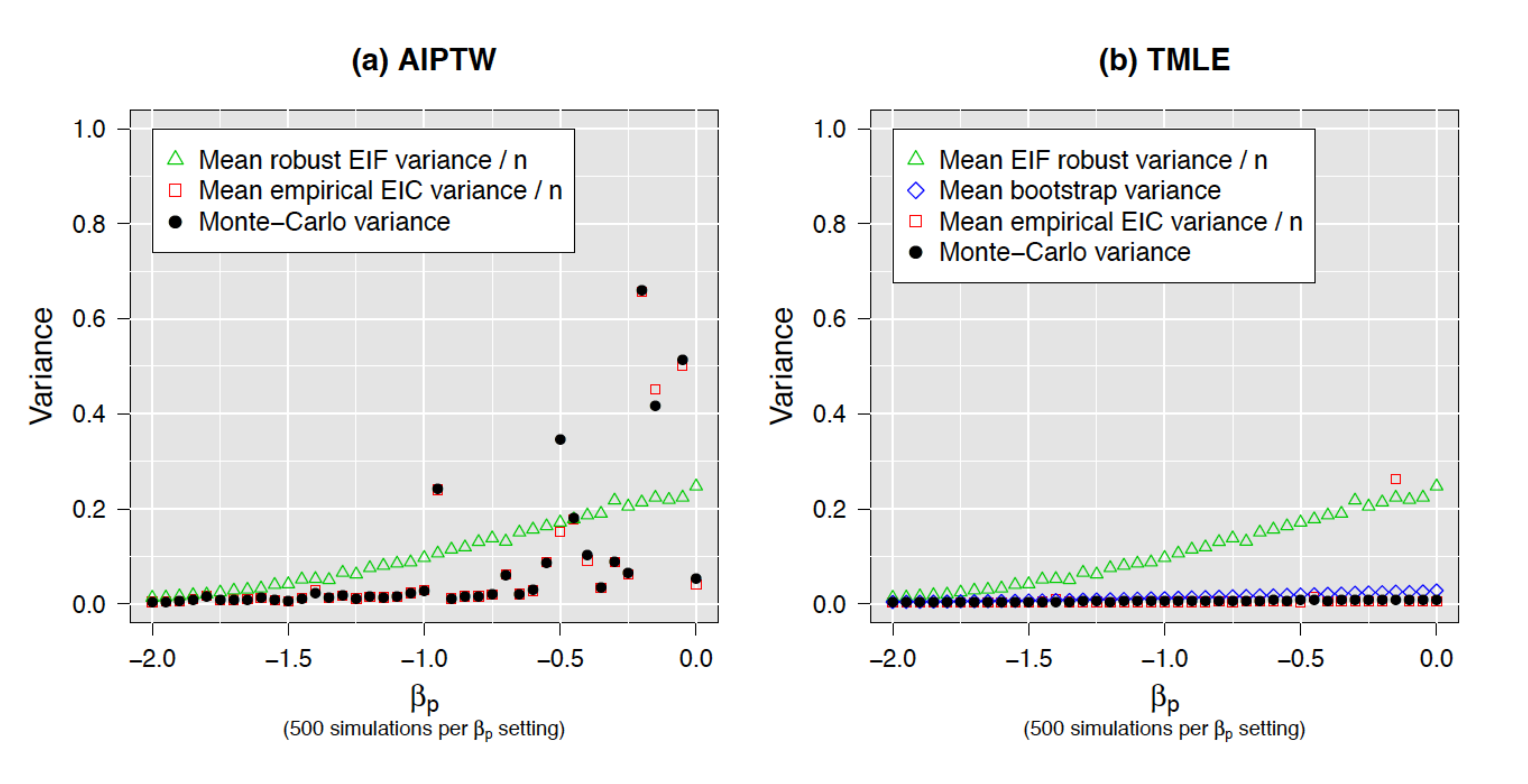}
\caption{Mean of variance estimates for each estimator under no treatment
effect ($\beta_{\psi_0}=0$) at each positivity ($\beta_{p}$) value under the
longitudinal treatment setting, overlaid with the estimator's Monte-Carlo
variance. Robust variance estimates are identical for the two estimators.}
\label{fig:compare2}
\end{figure}

Results for the longitudinal setting were less stable, though still similar to
the point treatment setting. Figure \ref{fig:compare2} shows the mean of the
variance estimates under each approach, overlaid with the Monte-Carlo variance
of the intervention specific mean outcome estimators. The same trend over the
different levels of positivity was seen as in Figure \ref{fig:compareVar}, with
the variance increasing with the magnitude of positivity issues. The empirical
EIF approach also performed well here at low levels of $\beta_{p}$ for both the
AIPW and TMLE estimators. At high values of $\beta_{p}$, the approach more
noticeably underestimate the variance of both intervention specific mean outcome
estimators. Consistent with the point treatment setting, the robust EIF approach
consistently over estimated the variance for both estimators. The bootstrap
approach resulted in slightly conservative variance, though were still very
similar to the Monte-Carlo variance estimates.

\begin{figure}[!h]
\centering
\includegraphics[width=5in]{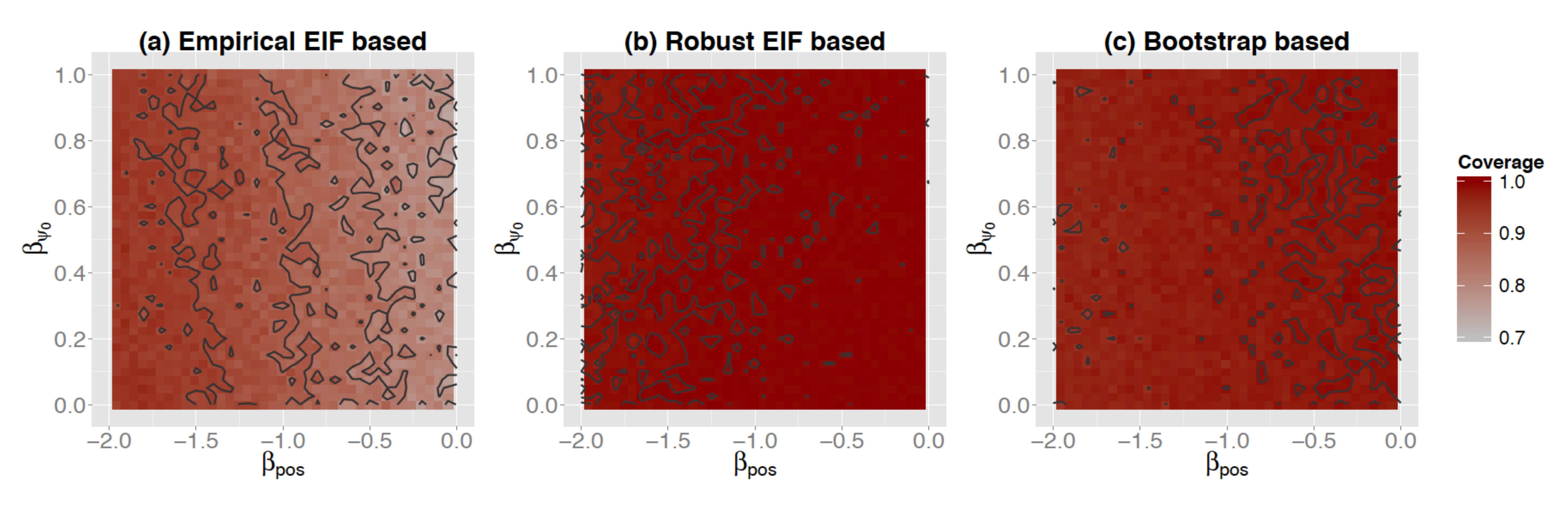}
\caption{Simulated coverage for each variance estimation approach for the TMLE
estimator under various treatment ($\beta_{\psi_0}$) and positivity
($\beta_{p}$) values under the longitudinal treatment setting.}
\label{fig:t3cover}
\end{figure}

Figure \ref{fig:t3cover} shows the coverage corresponding to each variance
estimation approach for the TMLE estimator of the intervention specific mean
outcome. Coverage for the empirical EIF approach dropped considerably with an
increase in positivity issues. For example, at a null effect (i.e.
$\beta_{\psi_0}$) the observed coverage was $0.93$ at $\beta_{p}=-2$ and $0.78$
at $\beta_{p}=0$. For the robust EIF approach, coverage increased with
positivity. This became as high as $1.00$ (i.e. all simulated confidence
intervals captured the true parameter value) at higher levels of positivity
issues. For the bootstrap approach, a higher level of coverage was also seen.
For example, under a null effect, a coverage of $0.95$ was observed at
$\beta_{p}=-2$ and $0.98$ at $\beta_{p}=0$.

Results for the Type I error and power were also similar to the point treatment
setting. When there was an effect, the empirical EIF approach resulted in the
highest power. At $\beta_{\psi_0}=1$ and $\beta_{p}=-2$, we observed a power of
$0.99$. However, the Type I error was also uncontrolled here, becoming as high
as $0.22$ at $\beta_{p}=0$. While the robust EIF approach maintained valid Type
I error rates, the power for this approach when an effect was present was the
lowest. For example, for an treatment effect size of $\beta_{\psi_0}=1$ we
observed a power ranging from $0.996$ at $\beta_{p}=-2$ to $0.14$ at
$\beta_{p}=0$. The bootstrap approach also resulted in controlled Type I error
rates, with observed values below $0.05$ over all values of $\beta_{p}$
considered. Power was higher than the robust EIF approach across all values of
$\beta_{\psi_0}$ and $\beta_{p}$. For a treatment effect size of
$\beta_{\psi_0}=1$, we observed a power ranging from $0.988$ at $\beta_{p}=-2$
to $0.40$ at $\beta_{p}=0$ for the bootstrap approach. Compared with the robust
EIF approach, this is almost a 3-fold increase in power.

\begin{figure}[!h]
\centering
\includegraphics[width=5in]{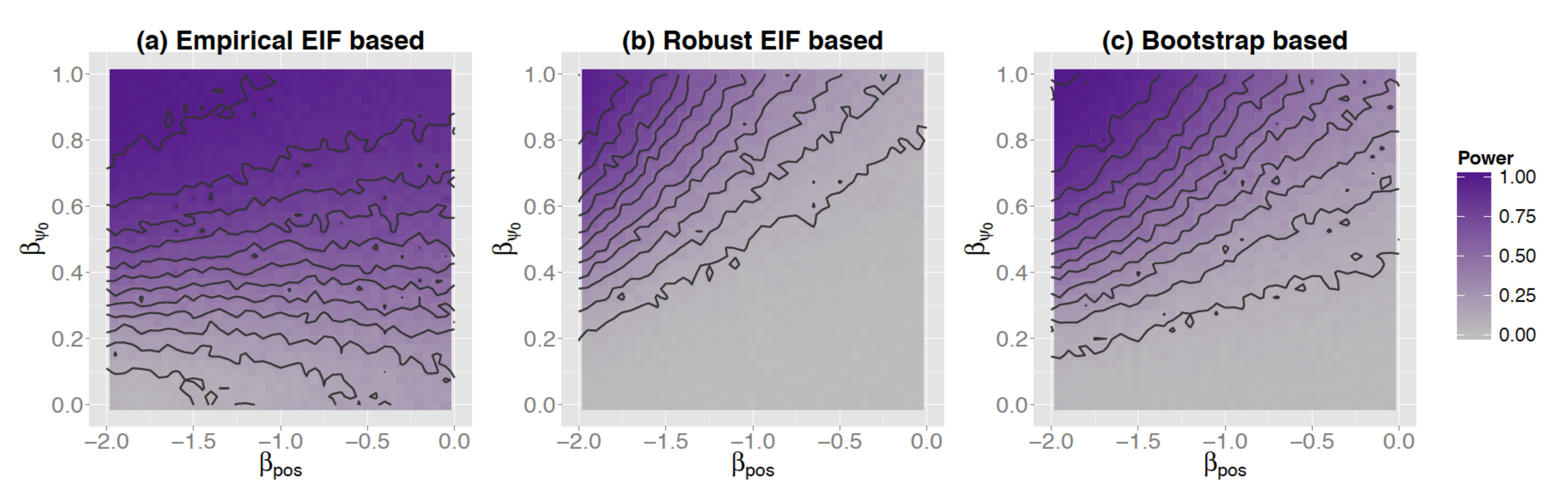}
\caption{Simulated power for each variance estimation approach for the TMLE
estimator under various treatment ($\beta_{\psi_0}$) and positivity
($\beta_{p}$) values under the longitudinal treatment setting.}
\label{fig:t3power}
\end{figure}

\section{Discussion}
\label{sec:disc}

The goal of the current study was to establish a consistent and robust approach
of estimating the variance of asymptotically efficient estimators such as TMLE,
estimating equations, and one step estimators which, in contrast to the common
approach based on the empirical variance of the estimated EIF, do not act
anti-conservatively when confronted with positivity violations. We have presented two
such approaches at estimating this variance: 1) a robust approach that directly
targets the asymptotic variance of the EIF, and 2) a bootstrap approach based on
fitting the initial density of the data once, followed by a non-parametric
bootstrap of the targeting step. In simulations, the variance of AIPW increases
noticeably with the variance of the EIF as positivity increases. The variance of
the TMLE was constrained in the face of increasing positivity violations, and as
a result, while the empirical EIF approach underestimated variance, the robust
EIF approach increasingly over-estimated the variance as the degree of
positivity violations increased. In contrast, the bootstrap based approach
provided less conservative variance estimation, while maintaining valid Type I
error control in the face of extreme positivity violations, both in the point
treatment and longitudinal setting.

We emphasize that, as the robust EIF approach directly targets the asymptotic
variance of the EIF, extremely large values of estimates from this can be used
to raise a red flag for unreliable statistical inference due to sparsity,
thereby declaring that the target parameter is practically not identifiable from
the data, and that the reported variance estimates (though large) will
themselves be imprecise. As such, we recommend that this approach be used if
there is concern regarding identifiability of the data for the target parameter
of interest.

While the EIF can raise a red flag for lack of identifiability, for substitution
estimators such as TMLE we suggest that it is overly conservative for
constructing valid confidence intervals and tests in finite sample in the face
of substantial positivity violations. In previous work \citep{Petersen2012}, we
suggested the parametric bootstrap as a diagnostic tool for sparsity-bias in the
point treatment setting. The approach can become cumbersome, as the analyst
would need to refit the whole likelihood for each iteration of the bootstrap.
Our robust EIF approach is able to avoid estimating the whole likelihood by
targeting the required means under the post intervention distribution defined by
the longitudinal g-computation formula directly. Even if we use an actual TMLE
of $P_0^d$, our analytic estimate of the variance is still much less computer
intensive than the parametric bootstrap method, in particular, in view that one
would need to run many replicate samples of the data set in order to pick up the
observations that correspond with a rare treatment and thus contribute large
numbers to the variance expression. Therefore, we believe that the proposed
analytic method will be (at least, practically) superior to the earlier proposed
parametric bootstrap method. Our presented bootstrap approach, while more
computationally intensive than the robust EIF approach, is also superior to the
earlier proposed approach, in that we do not have to refit the entire likelihood
under each iteration. This also significantly reduces the computational
resources required to obtain estimates of the target parameter, particularly if
computationally intensive non-parametric machine learning algorithms are used to
estimate these densities.

Further refinements can be applied in an attempt to obtain variance estimates
with an even smaller finite sample bias. One such approach is a convex
combination of the variance estimators considered above. For example, we noticed
in supplementary analyses that taking
\[\hat{\alpha}_n\hat{\sigma}^2_{eEIF,n}+(1-\hat{\alpha}_n)\hat{\sigma}^2_{rEIF,n}\]
had good performance, where $\hat{\sigma}^2_{eEIF,n}$ is the variance estimate
using the empirical EIF approach, $\hat{\sigma}^2_{rEIF,n}$ is the variance
estimate using the robust EIF approach, and
$\hat{\alpha}_n=|\hat{\sigma}^2_{rEIF,n}-\hat{\sigma}^2_{eEIF,n}| /
(\hat{\sigma}^2_{rEIF,n}+\hat{\sigma}^2_{eEIF,n})$. We note, however, that such
an approach is somewhat ad-hoc and may lead to varying results in other
simulations or distributions. We therefore chose not to present the results here.

A potential limitation of the robust approach at estimating the variance
involves the conditions for asymptotic linearity to be met. Note, however, that
we always require that our estimator (of the parameter) be asymptotically
linear. Thus, this is actually a limitation of our parameter estimator. Given
that we have an asymptotically linear estimator, we want a good estimator of its
variance. Furthermore, it is also required that $\bar{Q}^{d,\sigma^{2}_{t}}_{0}$
be estimated both consistently and at a fast enough rate. We limited the
computational complexity in our simulations by using simpler parametric models
to estimate $\bar{Q}^{d,\sigma^{2}_{t}}_{0}$, though a more non-parametric
approach such as Super Learning could have been applied. This approach can
become computationally expensive if there are many time points. In this regard,
the bootstrap approach is superior as it does not require the additional
estimation of $\bar{Q}^{d,\sigma^{2}_{t}}_{0}$.

It would be of interest to further evaluate not only the practical performance
of these variance estimation approaches in future studies, but also the
application of the approaches to other parameters. The appendix derives the
general approach for working marginal structural models. Further research into
the practical performance of this approach is needed for this setting. These
variance estimation approaches can also be used for more complex parameters,
such as mean outcomes under dynamic regimes, stochastic interventions, or
marginal structural working models. Future research could also develop a
collaborative TMLE \citep{VanderLaan2010} or cross-validated \citep{Zheng2010}
based approach at robustly estimating the EIF variance.

\clearpage
\bibliography{/Users/tranlm/Dropbox/Bibtex/school}

\clearpage
\section*{Appendix}

\subsection*{TMLE of $\sigma_{K+1}^2$ for marginal structural working models}
\label{sec:Dkplus1}
For the general working marginal structural model (MSM) $\Theta \equiv
\{m_{\beta}:\beta\}$ from Petersen el al. \citep{Petersen2014c}, we have that the component
corresponding with the last time point $K+1$ equals
\begin{eqnarray*}
D^*_{K+1}(P)&=&\sum_{d\in {\cal D}} h_1(d,K+1)
\frac{\mathbb{I}(\bar{A}(K)=d(\bar{L}(K)))
}{g_{0:K}(\bar{A}(K),\bar{L}(K))} (Y-\bar{Q}_{K+1}(\bar{A}(K),\bar{L}(K)))\\
&=& C_{K+1}(P)(\bar{A},\bar{L}) (Y-\bar{Q}_{K+1}),
\end{eqnarray*}
where, for some user defined weight function $h(d,K+1)$,
\begin{eqnarray*}
C_{K+1}(P)(\bar{A},\bar{L})&=&\sum_{d\in {\cal D}}
h_1(d,K+1)\frac{\mathbb{I}(\bar{A}(K)=d(\bar{L}(K))
)}{g_{0:K}(\bar{A},\bar{L})}. \\
h_1(d,K+1)&=&h(d,K+1)\frac{\frac{\partial}{\partial\beta}m_{\beta}(d,K+1)}{m_{\beta}(1-m_{\beta})}
\end{eqnarray*}

We want to obtain a representation of the variance of this component $D^*_{K+1}$
so that we can use a semi-substitution estimator of this part of the variance of
the EIF, hopefully, thereby obtaining a variance estimator
that is more accurate under violations of practical positivity, and a variance
estimator that can be used as a red flag for lack of identifiability.
This variance can thus be written as
\begin{equation*}
\begin{split}
\sigma^2_{K+1} &= \mathbb{E}[C^2(Y-\bar{Q}_{K+1})^2]\\
&=\mathbb{E} [C^2 \bar{Q}_{K+1}(1-\bar{Q}_{K+1})]\\
&=\mathbb{E}\left[ \left( \sum_{d \in {\cal D}} h_1(d,K+1)
\mathbb{I}(\bar{A}=d(\bar{L}))\right)^2
\frac{\bar{Q}_{K+1}(1-\bar{Q}_{K+1})}{g_{0:K}^2}(O) \right]\\
&=\mathbb{E}\left[ \left( \sum_{d_1,d_2 \in {\cal D}}
h_1(d_1,K+1)h_1(d_2,K+1)
\mathbb{I}(\bar{A}=d_1(\bar{L}))\mathbb{I}(\bar{A}=d_2(\bar{L}))\right)\frac{\bar{Q}(1-\bar{Q})}{g_{0:K}^2}(O) \right]\\
&=\sum_{d_1,d_2} h_1(d_1,K+1)h_1(d_2,K+1)
\mathbb{E}\left[\mathbb{I}(\bar{A}=d_1(\bar{L}))\mathbb{I}(\bar{A}=d_2(\bar{L}))\frac{
\bar{Q}(1-\bar{Q})}{g_{0:K}^2}(O)\right].
\end{split}
\end{equation*}
The latter expectation equals:
\[
\begin{array}{l}
\int_{\bar{L}} \mathbb{I}(d_1(\bar{L})=d_2(\bar{L})) \prod_{t=0}^{K+1}
q(L_t\mid \bar{A}(t^{-})=d_{1}(\bar{L}(t^{-})),\bar{L}(t^{-}))
\frac{\bar{Q}(1-\bar{Q})}{g_{0:t}(d_{1}(\bar{L}),\bar{L})}\\
= \mathbb{E}_{P_0^{d_1}}
\left[\mathbb{I}(d_1(\bar{L}_{d_1})=d_2(\bar{L}_{d_1}))
\frac{\bar{Q}(1-\bar{Q})}{g_{0:K}}(d_1(\bar{L}_{d_1}),\bar{L}_{d_1})\right].
\end{array}
\]
This yields the following expression:
\begin{align*}
\sigma^2_{K+1}&=\sum_{d_1,d_2 \in {\cal D}}
h_1(d_1,K+1)h_1(d_2,K+1) \\
&\hspace{.75in}\mathbb{E}\left[\mathbb{I}(d_1(\bar{L}_{d_1})=d_2(\bar{L}_{d_1}))\frac{
\bar{Q}(1-\bar{Q})}{g_{0:K}} (d_1(\bar{L}_{d_1}),\bar{L}_{d_1})\right] \\
&=\sum_{d_1 \in {\cal D}} h_1(d_1,K+1) \\
&\hspace{.75in} \mathbb{E}\left[ \left(
\sum_{d_2 \in {\cal D}}
h_1(d_2,K+1)\mathbb{I}(d_1(\bar{L}_{d_1})=d_2(\bar{L}_{d_1}))
\right) \frac{\bar{Q}(1-\bar{Q})}{g_{0:K}}
(d_1(\bar{L}_{d_1}),\bar{L}_{d_1})\right] \\
&=\sum_{d_1 \in {\cal D}} h_1(d_1,K+1) \mathbb{E}
Z_{d_1}(d_1,K+1) \numberthis \label{exprmsm}
\end{align*}
where
\[ Z(d_1,K+1) = \left(\sum_{d_2 \in {\cal D}} h_1(d_2,K+1)
\mathbb{I}(d_1(\bar{L}(K))=d_2(\bar{L}(K)))\right)
\frac{\bar{Q}(1-\bar{Q})}{g_{0:K}(d(\bar{L}(K)),\bar{L}(K))},\] so that the
counterfactual of $Z(d_1,K+1)$ under intervention $d_1$ is given by
\begin{align*}
Z_{d_1}(d_1,K+1)&=\left( \sum_{d_2 \in {\cal
D}}h_1(d_2,K+1)\mathbb{I}(d_1(\bar{L}_{d_1}(K))=d_2(\bar{L}_{d_1}(K)))
\right) \\
&\hspace{.75in} \frac{\bar{Q}(1-\bar{Q})}{g_{0:K}}
(d_1(L_{d_1}(K)),L_{d_1}(K)).
\end{align*}

\subsubsection*{Static regimens}

In the special case that the class of dynamic regimens ${\cal D}$ consists only
of static regimens $\bar{a}(K)$ so that there is only one and exactly one
treatment such that $\bar{A}(K)=d(\bar{L}(K))$, then we have \[ Z(K+1) =
h_1(\bar{A},K+1) \frac{\bar{Q}(1-\bar{Q})}{g_{0:K}}(\bar{A},\bar{L}),\]
so that
\[ Z_{d}(K+1)=h_1(d,K+1)\frac{\bar{Q}(1-\bar{Q})}{g_{0:K}}
(d(\bar{L}_d),\bar{L}_d).\] In that case, we have
\begin{eqnarray*}
\sigma^2_{K+1}&=& \sum_{d \in {\cal D}}
h_1(d,K+1)^2 \mathbb{E} Z_{1d}(K+1)
\end{eqnarray*}
where $Z_{1}(K+1)=
\bar{Q}(1-\bar{Q})/g_{0:K}(\bar{A},\bar{L})$ and $Z_{1d}(K+1)=
\bar{Q}(1-\bar{Q})/g_{0:K}(d(\bar{L}_d),\bar{L}_d)$.

It is important to note that in expressing our variance this way, we integrate
out the indicator of treatment over $\bar{A}$, i.e.
$\mathbb{I}(\bar{A}=d(\bar{L}))$. By getting rid of this indicator, we no longer
rely as heavily on observations from subjects following treatment in estimating
the variance of $D^{*}_{K+1}$. This particularly helps us when there is a lack
of positivity, since subjects with low probabilities of desired treatment simply
are not observed.
\newline

We have now shown that
\[
\sigma^2_{K+1}=\sum_{d \in {\cal D}} h_1(d,K+1)\mathbb{E}Z_{d}(d,K+1),\]
where
\[
Z(d_1,K+1)=\left\{ \sum_{d_2} h_1(d_2,K+1)
\mathbb{I}(d_1(\bar{L})=d_2(\bar{L})) \right\}
\frac{\bar{Q}(1-\bar{Q})}{g_{0:K}}(d_1(\bar{L}),\bar{L}).\]

We can now define $Z(K+1)(\bar{A},\bar{L})= \sum_{d \in {\cal D}}
h_1(d,K+1)\mathbb{I}(\bar{A}=d(\bar{L})) \frac{\bar{Q}(1-\bar{Q})}{g_{0:K}}
(\bar{A},\bar{L})$ (as function of $\bar{A},\bar{L}$) as a new outcome for our
longitudinal data structure such that
$Z_{d}(d,K+1)=Z(K+1)(d(\bar{L}_d),\bar{L}_d)$.
Our variance $\sigma^2_{K+1}$ is then represented as $\sum_{d \in {\cal D}}
h_1(d,K+1) \mathbb{E}Z_{d}(d,K+1)$. Thus, if we redefine the longitudinal data
as $(\bar{A},\bar{L})$ with the final outcome of interest as
$Z(K+1)=Z(K+1)(\bar{A},\bar{L})$, and use the working MSM parameter
$\mathbb{E}Z_{d}(K+1)=\beta_0$ with $\beta_0=\arg\min_{\beta} \sum_{d \in {\cal
D}} h_1(d,K+1) (\mathbb{E}Z_{d}(K+1)-\beta)^2$, then we have that
\[
\beta_0=\sum_{d \in {\cal D}} h_1(d,K+1) \mathbb{E}Z_{d}(K+1)/ \sum_{d \in {\cal
D}} h_1(d,K+1).
\]
This demonstrates that we can obtain $\sigma^2_{K+1}$ by simply multiplying
$\beta_0$ by $\sum_{d \in {\cal D}} h_1(d,K+1)$, i.e. \[ \sigma^2_{K+1}=\beta_0 \sum_d
h_1(d,K+1).
\] We can therefore also estimate this variance component $\sigma^2_{K+1}$ by
computing the TMLE of the intercept $\beta_0$ in the working MSM for our newly
defined outcome $Z(K+1)$ using weights $h_1(d,K+1)$, and then multiplying it
against $\sum_{d \in {\cal D}} h_1(d,K+1)$.

\subsection*{TMLE of $\sigma_{t}^2$ for marginal structural working models}

We now present the how to obtain a TMLE of the variance of the $t$-th component
of the EIF, $\sigma_{t}^2$. For the general working MSM from Petersen el al.
\citep{Petersen2014c}, we have that the component corresponding with the $t$-th
time point equals
\begin{eqnarray*}
D^*_{t}(P)&=&\sum_{d\in {\cal D}} h_1(d,t)
\frac{\mathbb{I}(\bar{A}(t^{-})=d(\bar{L}(t^{-})))}
{g_{0:t^{-}}(\bar{A}(t^{-}),\bar{L}(t^{-}))
}(\bar{Q}_{t^{+}}^d(\bar{A}(t),\bar{L}(t))-\bar{Q}_{t}^d)(\bar{A}(t^{-}),\bar{L}(t^{-})))\\
&=& \sum_{d\in {\cal D}} C_t(P,d) (\bar{Q}_{t^{+}}^d-\bar{Q}_t^d).
\end{eqnarray*}

Similar to above, we want to obtain a representation of the variance of this
component so that we can use a semi-substitution estimator of this part of the
variance of the EIF, hopefully, thereby obtaining a variance estimator that is
more accurate under violations of practical positivity, and a variance estimator
that can be used as a red flag for lack of identifiability.
This variance $\sigma^2_{t}$ can thus be written as
\begin{eqnarray*}
\sigma^2_t
&=&\sum_{d_1,d_2} h_1(d_1,t)h_1(d_2,t) \mathbb{E}\left[
\mathbb{I}(\bar{A}(t^{-})=d_1)\mathbb{I}(\bar{A}(t^{-})=d_2)
\frac{\Sigma_t(d_1,d_2)}{g^2_{0:t^{-}}}(\bar{A}(t^{-}),\bar{L}(t^{-}))\right]
\end{eqnarray*}
where
\[
\Sigma_t(d_1,d_2)(\bar{A}(t^{-}),\bar{L}(t^{-}))=\mathbb{E}\left [
(\bar{Q}^{d_1}_{t^{+}}-\bar{Q}^{d_1}_{t})(\bar{Q}^{d_2}_{t^{+}}-\bar{Q}^{d_2}_{t})\mathrel{\big|}
\bar{A}(t^{-}),\bar{L}(t^{-})\right]
\]
is the conditional covariance of $\bar{Q}^{d_1}_{t^{+}}$ and
$\bar{Q}^{d_2}_{t^{+}}$, given $(\bar{A}(t^{-}),\bar{L}(t^{-}))$. Note that this
can be obtained by regressing this cross-product on
$(\bar{A}(t^{-}),\bar{L}(t^{-}))$.
The latter sum can be further worked out giving us
\begin{eqnarray*}
\sigma^2_t&=& \sum_{d_1 \in {\cal D}} h_1(d_1,t) \mathbb{E}Z_{d_1}(d_1,t),
\end{eqnarray*}
where
\[
Z(d_1,t)=\left(\sum_{d_2 \in {\cal D}}h_1(d_2,t)
\mathbb{I}(d_1(\bar{L}(t^{-}))=d_2(\bar{L}(t^{-})) )\right)\frac{\Sigma_t(d_1,d_2) }{
g_{0:t^{-}}}(d_{1,t^{-}}(\bar{L}(t^{-})),\bar{L}(t^{-})).\] so that the
counterfactual of $Z_t$ under intervention $d_1$ is given by
\[
Z_{d_1}(d_1,t)= \left( \sum_{d_2 \in {\cal D}}h_1(d_2,t)
\mathbb{I}(d_1(\bar{L}_{d_1}(t^{-}))=d_2(\bar{L}_{d_1}(t^{-})))
\right)\frac{\Sigma_t(d_1,d_2)}
{g_{0:t^{-}}}(d_{1,t^{-}}(\bar{L}_{d_1}(t^{-})),\bar{L}_{d_1}(t^{-})).
\]

With this expression, we can now use a TMLE to estimate
$\mathbb{E}Z_{d_1}(d_1,t)$ for each $d_1\in {\cal D}$ by using the longitudinal
data structure with final outcome $Z(d_1,t)$, for each $d_1$ separately. To
create the observed outcome $Z(d_1,t)$ we need a fit of the treatment mechanism
$g_{A(m)}:m=0,1,\ldots,t^{-}$, evaluated at
$\bar{A}(t^{-})=d_{t^{-}}(\bar{L}(t^{-}))$, and for each rule compatible with
$d_1$ (for that subject) we need to have an estimate of $\Sigma_t(d_1,d_2)$.
Thus, given a priori estimates of the full treatment mechanism and all
$(\Sigma_t(d_1,d_2):d_1,d_2 \in {\cal D})$ we can construct this observed
outcome $Z(d_1,t)$ and run the TMLE.

\subsection*{Estimation of the variance of the EIF}

The above approach defines for each time point $t$ and each rule $d$ an observed
longitudinal outcome $Z(d,t)$, where $Z(d,t)$ is a function of
$(\bar{A}(t),\bar{L}(t))$. The TMLE of $\mathbb{E}Z_{d}(d,t)$ can then be
computed based on the longitudinal data structure
$(L(0),A(0),\ldots,L(t),A(t),Z(d,t))$ for each $d$ and each $t \in \{0, 1,
\ldots,K+1\}$. As a result, we have that $\sigma^2_t=\sum_{d \in {\cal D}}
h_1(d,t)\mathbb{E}Z_{d}(d,t)$ and
\begin{eqnarray*}
\sigma^2&=&\sum_{t=0}^{K+1} \sigma^2_t \\
&=&\sum_{d \in {\cal
D}} \left( \sum_{t=0}^{K+1} h_1(d,t)\mathbb{E}Z_{d}(d,t)\right) \\
&=&\sum_{d \in {\cal D}} \mathbb{E}\left[ \sum_{t=0}^{K+1}
h_1(d,t)Z_{d}(d,t)\right].
\end{eqnarray*}
Let's now define the counterfactual outcome
\[\bar{Z}_d(d)\equiv \sum_{t=0}^{K+1} h_1(d,t)Z_{d}(d,t),\]
and the corresponding observed outcome
\[\bar{Z}(d)\equiv \sum_{t=0}^{K+1} h_1(d,t)Z(d,t).\]
We could apply the TMLE to estimate $\mathbb{E}\bar{Z}_d(d)$  based on the
longitudinal data structure $(L(0),A(0),\ldots,L(K),A(K),\bar{Z}(d,K+1))$, for
each $d\in {\cal D}$, and use that \[\sigma^2=\sum_{d \in {\cal D}}
\mathbb{E}\bar{Z}_d(d).\]

In applying TMLE here, we should be using that \[
\mathbb{E}\left[\bar{Z}_d \mathrel{\big|} \bar{A}(m),\bar{L}(m)\right]=\sum_{t\leq
m}h_1(d,t)Z(d,t)+ \mathbb{E}\left[\sum_{t>m}h_1(d,t)Z(d,t) \mathrel{\Big|} \bar{A}(m),\bar{L}(m)\right].\]
To start with, let
\begin{eqnarray*}
\bar{Q}_d^{Z(K+1)}&=&\mathbb{E}[\bar{Z}(d) \mid \bar{A}(K),\bar{L}(K)] \\
&=&\sum_{t \leq K}h_1(d,t)Z(d,t)+\mathbb{E}[h_1(d,K+1)+Z(d,K+1)\mid
\bar{A}(K),\bar{L}(K)].
\end{eqnarray*}
Denote last conditional expectation with $\bar{Q}_d^{Z(K+1),d}$ so that
\[
\bar{Q}_d^{Z(K+1)}=\sum_{t\leq K}h_1(d,t)Z(d,t)+\bar{Q}_d^{Z(K+1),d}.\]
Then,
\begin{eqnarray*}
\bar{Q}_d^{Z(K)}&=&\mathbb{E}\left[\bar{Q}_d^{Z(K+1)}\mathrel{\Big|}
\bar{A}(K-1),\bar{L}(K-1)\right]\\
&=&\sum_{t\leq K-1}h_1(d,t)Z(d,t)+
\mathbb{E}\left[h_1(d,K)Z(d,K)+\bar{Q}_d^{Z(K+1),d} \mathrel{\Big|}
\bar{A}(K-1),\bar{L}(K-1)\right].
\end{eqnarray*}
Again, denote the latter conditional expectation by $\bar{Q}_d^{Z(K),d}$ so that
\[\bar{Q}_d^{Z(K)}=\sum_{t\leq K-1}h_1(d,t)Z(d,t)+\bar{Q}_d^{Z(K),d}.\]
Then,
\begin{eqnarray*}
\bar{Q}_d^{Z(K-1)}&=&\mathbb{E}\left[\bar{Q}_d^{Z(K)} \mathrel{\Big|}
\bar{A}(K-2),\bar{L}(K-2)\right] \\
&=&\sum_{t\leq K-2}h_1(d,t)Z(d,t)+ \\
&&\hspace{.5in}\mathbb{E}\left[h_1(d,K-1)Z(d,K-1)+\bar{Q}_d^{Z(K),d}
\mathrel{\Big|} \bar{A}(K-2),\bar{L}(K-2)\right].
\end{eqnarray*}
Again, denote the latter conditional expectation with $\bar{Q}_d^{Z(K-1),d}$ so
that
\[
\bar{Q}_d^{Z(K-1)}=\sum_{t\leq K-2}h_1(d,t)Z(d,t)+\bar{Q}_d^{Z(K-1),d}.\]
This is then iterated:
\[
\bar{Q}_d^{Z(m)}=\sum_{t\leq m-1}h_1(d,m)Z(d,m)+\bar{Q}_d^{Z(m),d},\]
where
$\bar{Q}_d^{Z(m),d}=\mathbb{E}\left[h_1(d,m)Z(d,m)+\bar{Q}_d^{Z(m+1),d}
\mathrel{\Big|} \bar{A}(m-1),\bar{L}(m-1)\right]$.

Before we go to the next conditional expectation we need to target with a
parametric submodel, such as
\[ \mbox{Logit}\,\bar{Q}_d^m(\epsilon)=\mbox{Logit}\,\bar{Q}_d^m+\epsilon
\frac{\mathbb{I}(\bar{A}(m-1)=d(\bar{L}(m-1)))}{g_{0:m-1}}.\] In this way, we
will only have to run one TMLE for each rule $d$, which still utilizes that the
outcome is a sum of outcomes that are known for histories including that
outcome.

\end{document}